\magnification=\magstep1
\input amstex
\documentstyle{amsppt}
\pageheight{9truein}
\loadbold
\loadmsbm
\loadeufm
\UseAMSsymbols
\baselineskip=12pt
\def\var{\varepsilon}
\def\bC{\Bbb C}
\def\bN{\Bbb N}
\def\bR{\Bbb R}
\def\bS{\Bbb S}
\def\bZ{\Bbb Z}
\def\cO{\Cal O}
\def\bO{\Bbb O}
\def\cA{\Cal A}
\def\cE{\Cal E}
\def\cH{\Cal H}
\def\cF{\Cal F}
\def\cJ{\Cal J}
\def\cB{\Cal B}
\def\cS{\Cal S}
\def\cC{\Cal C}
\def\cK{\Cal K}
\def\oD{\overline D}
\def\bP{\Bbb P}
\def\id{\text{id}}
\def\Im{\text{Im }}
\def\diam{\text{ diam}}
\def\Ker{\text{Ker }}
\def\Coker{\text{Coker }}
\def\Hom{\text{Hom}}
\def\bHom{\bold{Hom}}

\def\Ob{\text{Ob }}
\def\pr{\text{pr}}
\def\supp{\text{ supp}}
\NoBlackBoxes
\topmatter
\title Analytic Sheaves in Banach Spaces\endtitle
\address Department of Mathematics
Purdue University
West Lafayette, IN\ 47907-1395\endaddress
\author
L\'aszl\'o Lempert\endauthor\footnote""{Research partially supported by NSF
grant DMS 0203072\hfill\break}
\abstract
We introduce a class of analytic sheaves in a Banach space $X$, that we shall
call cohesive sheaves.
Cohesion is meant to generalize the notion of coherence from finite dimensional
analysis.
Accordingly, we prove the analog of Cartan's Theorems A and B for cohesive
sheaves on pseudoconvex open subsets $\Omega\subset X$, provided $X$ has an
unconditional basis.
\endabstract
\endtopmatter
\document
\TagsOnRight
\head Introduction\endhead

In finite dimensional complex analysis and geometry coherent analytic sheaves
play a central role, for the following four reasons:
\item\item{(i)} Most sheaves that occur in the subject are coherent.
\item\item{(ii)} Over pseudoconvex subsets of $\bC^n$ their higher
cohomology groups vanish.
\item\item{(iii)} The class of coherent sheaves is closed under natural
operations.
\item\item{(iv)} Whether a sheaf is coherent can be decided locally.

The purpose of this paper is to introduce a comparable class of sheaves, that we
shall call cohesive, in Banach spaces.
This notion is different from coherence, which formally makes sense in infinite
dimensions as well.
However, coherence is not relevant for infinite dimensional geometry, since it
has to do with finite generation, while in infinite dimensional spaces one
frequently encounters sheaves, such as tangent sheaves and ideal sheaves of
points, that are not finitely generated over the structure sheaf.
The structure sheaf itself is not known to be coherent in any infinite
dimensional Banach space, either.

We shall define cohesive sheaves 
in general Banach spaces (always over $\bC$).
However, we are able to prove meaningful results about cohesive sheaves only in
some Banach spaces, e.g.~in those that have an unconditional basis.
(Our main results do not generalize to certain nonseparable spaces, and we do
not know whether they hold in all separable spaces, or at least in those that
have a Schauder basis.---For the notion of Schauder and unconditional bases, see Section~1.)

Cohesive sheaves are sheaves of modules with an extra structure and a special
property.
We shall arrive at their definition in four steps.
Given a Banach space $X$, an open $\Omega\subset X$, and another Banach 
space $E$,
the sheaf $\cO^E=\cO^E_\Omega$ of germs of holomorphic functions 
$(\Omega,x)\to E$, $x\in\Omega$, will be called a plain sheaf.
It is regarded as a sheaf of modules over $\cO=\cO^\bC$.
If $U\subset\Omega$ is open, $F$ is another Banach space and 
$\Hom(E,F)$ denotes the Banach space of bounded 
linear
operators, then any holomorphic $\varphi\colon U\to\Hom(E,F)$ 
induces a homomorphism $\cO^E|U\to\cO^F|U$.
Such homomorphisms will be called plain.

Next one defines when a sheaf $\cA$ of $\cO$--modules over $\Omega$ is analytic.
In traditional finite dimensional terminology analytic sheaves are the same as
sheaves of $\cO$--modules.
However, the latter notion is adequate only if one is satisfied with studying
finitely generated $\cO$--modules; already in finite dimensional spaces the
sheaf of Banach (space or bundle) valued holomorphic germs has a richer
structure than a mere $\cO$--module. For example, certain infinite sums
of sections make sense, which cannot be explained in terms of the
$\cO$--module structure.
In the finite dimensional context Leiterer in [Li] proposed to capture this
richness by introducing a Fr\'echet space structure on $\Gamma(U,\cA)$,
$U\subset \Omega$ open; in this way he obtained the notion of an analytic Fr\'echet sheaf.
For infinite dimensional $\Omega$ the corresponding notion would not be
practical, though, and instead our definition will be inspired by a suggestion
of Douady [D1].
To generalize the notion of a complex analytic space, he proposed that the
infinite dimensional analog of a ringed space should be a ``functored space''.
We shall say that a sheaf $\cA$ of $\cO$--modules over $\Omega$ is endowed with
an analytic structure if for every plain sheaf $\cE$ a submodule {\bf
Hom}$(\cE,\cA)\subset\bHom_{\cO}(\cE,\cA)$ is specified.
The correspondence $\cE\mapsto\bHom(\cE,\cA)$ should satisfy two natural conditions;
then we also say that $\cA$ is an analytic sheaf.
An $\cO$--homomorphism $\cE\to\cA$ will be called analytic if its
germs are in $\bHom(\cE,\cA)$.
Any plain sheaf $\cA$ has a natural analytic structure, namely $\bHom(\cE,\cA)$ will consist of germs of plain homomorphisms.

Our notion of analyticity slightly conflicts with the traditional terminology:
while every sheaf of $\cO$--modules admits an analytic structure, this structure is not unique, see 3.7.

Now consider an infinite sequence
$$
\ldots\to\cF_2\to\cF_1\to\cA\to 0\tag*
$$
of analytic sheaves and homomorphisms over $\Omega$, with each $\cF_n$ plain.
We shall say that (*) is a complete resolution of $\cA$ if for each pseudoconvex $U\subset\Omega$ the induced sequence on sections
$$
\ldots\to\Gamma (U,\cF_2)\to\Gamma(U,\cF_1)\to\Gamma(U,\cA)\to 0\tag**
$$
is exact and, moreover, the same is true if in (**) each $\cF_n$ is 
replaced by $\bHom(\cE,\cF_n)$ and $\cA$ by $\bHom(\cE,\cA)$, for any plain 
sheaf $\cE$.

Finally, we shall call an analytic sheaf $\cA$ over $\Omega$ cohesive if 
$\Omega$ can be covered by open sets over which $\cA$ has a complete 
resolution.

Analytic sheaves locally isomorphic to plain sheaves obviously are 
cohesive, but at first glance the notion of a complete resolution is 
quite formidable, and it is not clear how such a resolution can be 
constructed beyond trivial situations.
An exact sequence like (*) is perhaps doable, but how will one ensure the 
exactness of (**) {\sl before} cohomology vanishing is known?
An answer is that in the spaces we consider it suffices to check the 
exactness of (**) for $U$ in an appropriate subclass of pseudoconvex sets, 
see Section 6.
Using this device it is possible to show that certain important ideal sheaves 
are cohesive (see Section 10):

\proclaim{Theorem~1}Suppose a Banach space $X$ has an unconditional basis, 
$\Omega\subset X$ is open, and $M\subset \Omega$ is a direct submanifold 
(see Section~1).
If $E$ is another Banach space, then the sheaf $\cJ\subset\cO^E_{\Omega}$ 
of germs that vanish on $M$ is cohesive.
\endproclaim

The main result of this paper is the following generalization of Cartan's 
Theorems A and B (see Section 9):

\proclaim{Theorem~2}If a Banach space $X$ has an unconditional basis, 
$\Omega\subset X$ is pseudoconvex, and $\cA$ is a cohesive sheaf over 
$\Omega$, then
\itemitem{(a)}$\cA$ admits a complete resolution over $\Omega$; and
\itemitem{(b)}$H^q(\Omega,\cA)=0$ for $q\geq 1$.
\endproclaim

If coherence was useful because of the four reasons listed earlier, it would 
make sense to test the notion of cohesion against (i),...,(iv).
As to (iv) and (ii), being a cohesive sheaf is a local property, and at least in spaces with unconditional bases, higher cohomology groups of cohesive sheaves do vanish.
As to (i), only the future will tell exactly what sheaves will occur in 
infinite dimensional analysis and geometry, and whether they are cohesive.
However, again in spaces with unconditional bases, two clearly important types 
of sheaves are cohesive:\ the sheaf of sections of locally trivial holomorphic 
Banach bundles, and ideal sheaves of direct submanifolds. For this reason,
the theory of cohesive sheaves allows one to study geometry and analysis
on direct submanifolds.

In the Appendix we show that the ideal sheaf of certain analytic subsets is 
not cohesive; but, since those analytic subsets are  ``pathological'' 
anyway,  this result  counts in favor of the notion of 
cohesion. Indeed, it suggests that intuitively bizarre analytic subsets can 
be eliminated from complex geometry by testing whether their ideal sheaf is 
cohesive. We shall briefly introduce and study analytic subvarieties along 
these lines in Section~11.

Finally, item (iii) above.
The class of cohesive sheaves is closed under certain sheaf theoretical operations, but not under all operations that coherent sheaves admit; we give various examples as we develop the theory.
This seems to be a feature one has to learn to live with, and is due to peculiarities that Banach space valued functions can exhibit, even those of finitely many variables.

A few words about the history of analytic cohomology vanishing in vector spaces.
Cartan's theorems were generalized by Bishop and Bungart, and then by 
Leiterer 
to certain, so called Banach coherent analytic Fr\'echet sheaves over 
finite dimensional Stein spaces.
Douady proved the vanishing of higher cohomology groups of certain 
sheaves over compact subsets of Banach spaces; see [Bi,Bu,C,D1,Li].
More recently we considered pseudoconvex sets $\Omega$ in Banach spaces that 
have an unconditional basis, and proved in [L4] that for a trivial Banach 
(or even Fr\'echet) bundle over $\Omega$ higher cohomology groups vanish.
This was generalized in [P2-3], and finally in [L6] to arbitrary locally 
trivial holomorphic Banach bundles.
Further results on sheaf cohomology were obtained by Patyi in [P4-5].
The present paper borrows ideas from most of these works.

While putting the finishing touches on this paper, I received two related
preprints from Patyi. Among other things, [P6] features a cohomology 
vanishing
result for holomorphic Banach bundles, whose assumptions are probably
less restrictive than the assumptions of [L4] and of 9.1 of the
present paper.
In [P7] Patyi considers the sheaf $\cJ$ of Theorem 1 above, when 
$M=\Omega\cap Y$ and $Y\subset X$ is a complemented subspace. Under the
assumptions of [P6], he constructs a resolution of $\cJ$ that in essence
shows that $\cJ$ is cohesive---although he does not introduce this 
notion---, and
proves $H^q(\Omega,\cJ)=0$ for $q\ge 1$.

Acknowledgement.
I am grateful to J.~Lipman, J.~McClure, E. Szab\'o, and L. Tong for 
discussions we have had on various aspects of this work.

\head 1.\ Glossary\endhead
\subhead 1.1\endsubhead
For matters of sheaf theory we refer to [Se] or [Br], for complex analysis 
in Banach spaces to [M], and for basic notions of infinite dimensional 
complex geometry to [L1, Sections 1--2].
In this glossary we shall nevertheless spell out a few basic definitions of 
infinite dimensional complex analysis and geometry.
Let $X$ and $E$ be Banach spaces, always over $\bC$, and $\Omega\subset X$ 
open.
We denote the space of continuous linear operators $X\to E$ by $\Hom(X,E)$, 
endowed with the operator norm.

\definition{1.2.\ Definition}A function 
$f\colon\Omega\to E$ is holomorphic 
if for each $x\in\Omega$ there is an $L\in\Hom(X,E)$ such that
$$
f(y)-f(x)=L(y-x)+o(\|y-x\|_X),\qquad \Omega\ni y\to x.
$$
\enddefinition

This leads to the notion of what we call in [L1] a rectifiable complex Banach 
manifold:\ it is a Hausdorff space, sewn together from open subsets of Banach 
spaces (charts), with holomorphic sewing functions.
Using the charts one can define  holomorphic maps between rectifiable complex Banach manifolds.

\subhead 1.3\endsubhead
Let $M$ be a rectifiable complex Banach manifold.

\definition{Definition}(a) A closed subset $N\subset M$ is a submanifold if for each $x\in N$ there are a neighborhood $U\subset M$, an open subset $O$ of a Banach space $X$, a closed subspace $Y\subset X$, and a biholomorphic map $U\to O$ that maps $U\cap N$ onto $O\cap Y$.

(b) $N$ is a direct submanifold if, in addition, $Y$ above has a closed complement.
\enddefinition

The definition implies that the submanifold $N$ itself is a rectifiable complex Banach manifold.

\definition{Definition}A locally trivial holomorphic Banach bundle over $M$ is a holomorphic map $\pi\colon L\to M$, where $L$ is a rectifiable complex Banach manifold, and each fiber $\pi^{-1}x$, $x\in M$, is endowed with the structure of a topological vector space.
It is required that for every $x\in M$ there be a neighborhood $U\subset M$, 
a Banach space $E$, and a biholomorphic map $\pi^{-1}U\to E\times U$ that 
maps each fiber $\pi^{-1} y$, $y\in M$, isomorphically on 
$E\times\{y\}\approx E$.
\enddefinition

\definition{1.4.\ Definition}An upper semicontinuous function $u\colon\Omega\to [-\infty,\infty)$ is plurisubharmonic if for every pair $x,y\in X$ the function $\bC\ni\lambda\mapsto u(x+\lambda y)$ is subharmonic where defined.
\enddefinition

\definition{Definition}The open set $\Omega$ is pseudoconvex if for each 
finite dimensional subspace $Y\subset X$ the set $\Omega\cap Y$ is 
pseudoconvex in $Y$.
\enddefinition

It follows from the characterization of pseudoconvexity in
[M, 37.5 Theorem (e) or (f)] that if $\Omega$ is pseudoconvex, then so is any 
$\Omega'\subset X'$ biholomorphic to it.

\definition{1.5.\ Definition}In a Banach space $X$ a sequence $e_1,e_2,\ldots$ is a Schauder basis if every $x\in X$ can be uniquely represented as a norm convergent sum
$$
x=\sum^\infty_{n=1}\lambda_n e_n,\qquad \lambda_n\in\bC.
$$
The basis is unconditional if, in addition, the above series converge after arbitrary rearrangements.
\enddefinition

The spaces $l^p$, $L^p[0,1]$ for $1<p<\infty$, $l^1$, and the space $c_0$ of sequences converging to zero, all have unconditional bases, but $L^1[0,1]$ and $C[0,1]$ have Schauder bases only, see [Sn].

\head 2.\ Plain sheaves\endhead
\subhead 2.1\endsubhead
Fix a Banach space $(X,\|\ \|)$ and an open $\Omega\subset X$.
This notation will be used throughout the paper.
As explained in the Introduction, if $E$ is another Banach space, the sheaf $\cO^E=\cO^E_{\Omega}$ of germs of holomorphic functions 
$(\Omega,x)\to E$, $x\in\Omega$, is called a plain sheaf.
In particular, $\cO=\cO^{\bC}$ is a sheaf of rings and each $\cO^E$ is in a natural way a sheaf of $\cO$--modules.
In general, a sheaf of $\cO$--modules often will simply be called an $\cO$--module.
The sheaf of homomorphisms between $\cO$--modules $\cA$ and $\cB$ will be denoted $\bHom_{\cO}(\cA,\cB)$, itself an $\cO$--module.

If $F$ is also a Banach space and $U\subset \Omega$ is open, any 
holomorphic $\varphi\colon U\to\Hom(E,F)$ defines a homomorphism 
$\cO^E|U\ni\bold e\mapsto\varphi\bold e\in \cO^F|U$.
Such homomorphisms are called plain homomorphisms. 
The following is obvious:

\proclaim{Proposition}If a plain homomorphism induced 
by $\varphi\colon U\to\Hom(E,F)$ annihilates germs of 
constant functions, then $\varphi=0$.
\endproclaim

It follows that germs of plain homomorphisms $\cO^E|U\to\cO^F|U$
form an $\cO$--module
$$
\bHom_{\text{plain}}(\cO^E,\cO^F)\subset\bHom_\cO(\cO^E,\cO^F),
$$
isomorphic to $\cO^{\Hom(E,F)}$; its sections over any open $U\subset\Omega$
are in one--to--one correspondance with the plain homomorphisms 
$\cO^E|U\to\cO^F|U$.

\head 3.\ Analytic sheaves\endhead
\subhead 3.1\endsubhead If $\cA$ is an $\cO$--module over $\Omega$, an analytic structure on $\cA$ is the choice, for each plain sheaf $\cE$, of a submodule $\bHom(\cE,\cA)\subset\bHom_{\cO}(\cE,\cA)$, subject to

\itemitem{(i)}if $\cE$ and $\cF$ are plain sheaves, $x\in\Omega$, and 
$\boldsymbol\varphi\in\bHom_{\text{plain}}(\cE,\cF)_x$, then\newline
${\boldsymbol\varphi}^*\bHom(\cF,\cA)_x\subset\bHom(\cE,\cA)_x$;
\itemitem{(ii)}$\bHom(\cO,\cA)=\bHom_{\cO}(\cO,\cA)$.

An $\cO$--module endowed with an analytic structure is called an analytic sheaf.
The restriction of an analytic sheaf to an open $U\subset\Omega$ inherits an analytic structure.
Given analytic sheaves $\cA,\cB$ over $\Omega$ and an open $U\subset\Omega$, a homomorphism $\varphi\colon\cA|U\to\cB|U$ is called analytic if for every plain $\cE$ the induced homomorphism
$$
\varphi_*\colon\bHom_{\cO}(\cE|U,\cA|U)\to\bHom_{\cO}(\cE|U,\cB|U)
$$
maps $\bHom(\cE|U,\cA|U)$ in $\bHom(\cE|U,\cB|U)$.

Any plain sheaf $\cF$ can be endowed with an analytic structure 
$\bHom(\cE,\cF)=\bHom_{\text{plain}}(\cE,\cF)$.
In what follows, we shall automatically use this analytic structure on 
plain sheaves.
More generally, the sheaf $\cS$ of holomorphic sections of a locally trivial 
holomorphic Banach bundle $S\to\Omega$ has a canonical analytic structure:\ 
if over an open $U\subset\Omega$ the restriction $S|U$ is isomorphic to a 
trivial bundle $F\times U\to U$, this isomorphism induces an isomorphism of 
$\cO$--modules
$$
\bHom_{\cO}(\cO^E,\cO^F)|U\overset\sim\to\rightarrow \bHom_{\cO}(\cO^E,\cS)|U,
$$
and we define $\bHom(\cO^E,\cS)$ so that its restriction to such $U$ is the 
image of $\bHom(\cO^E,\cO^F)|U$.

Analytic homomorphisms between plain sheaves are just the plain homomorphisms, 
and, more generally, a homomorphisms between plain, resp.~analytic sheaves 
$\cE$ and $\cB$ is analytic if its germs are in $\bHom(\cE,\cB)$.
We shall write $\bHom(\cA,\cB)$ for the sheaf of germs of 
analytic homomorphisms between analytic sheaves $\cA$ and $\cB$; if $\cA$ is 
plain, this is consistent with notation already in use.

It is possible to define an analytic structure on the sheaf $\bHom(\cA,\cB)$ itself, but in this generality the structure will not have useful properties.
To keep the discussion simple we shall therefore not consider $\bHom(\cA,\cB)$ as an analytic sheaf.

\subhead 3.2\endsubhead
The notion of an analytic structure is naturally expressed in the language of 
category theory---even in more than one way---, of which we shall avail ourselves only very sparingly.
Consider the category $\bS$ of $\cO$--modules and $\cO$--homomorphisms over $\Omega$ (a monoidal category with unit $\cO$).
Now $\cO$--modules over $\Omega$ also form a category {\sl enriched} over 
$\bS$, or an $\bS$--category, see [K], that we shall denote $\bO$.
This means that with any pair $\cA,\cB\in\Ob\bO$, i.e.~$\cO$--modules, an $\cO$--module $\bHom_{\cO}(\cA,\cB)\in\bS$ is associated (the ``hom--object''); further, homomorphisms
$$
\bHom_{\cO}(\cA,\cB)\otimes\bHom_{\cO}(\cB,\cC)\to\bHom_{\cO}(\cA,\cC)
$$
are specified, satisfying certain axioms.
There is also the $\bS$--subcategory $\bP$ of $\bO$, whose objects are plain sheaves and the hom--object associated with $\cE,\cF\in\Ob\bP$ is 
$\bHom_{\text{plain}} (\cE,\cF)$.

Any $\cA\in\Ob\bO$ determines a contravariant $\bS$--functor $\bHom_{\cO}(\cdot,\cA)$ from $\bP$ to $\bO$.
This functor associates with $\cE,\cF\in\Ob\bP$ the homomorphism
$$
\bHom_{\text{plain}}(\cE,\cF)\to
\bHom_{\cO}(\bHom_{\cO}(\cF,\cA),\bHom_{\cO}(\cE,\cA)),
$$
induced by composition.
In this language an analytic structure on $\cA\in\Ob\bO$ is an $\bS$--subfunctor $\bHom(\cdot,\cA)$ of $\bHom_{\cO}(\cdot,\cA)$ satisfying $\bHom(\cO,\cA)=\bHom_{\cO}(\cO,\cA)$; subfunctor meaning that $\bHom(\cE,\cA)\subset\bHom_{\cO}(\cE,\cA)$ is a submodule for every $\cE\in\Ob\bP$.

\subhead 3.3\endsubhead
The sum $\cA\oplus\cB$ of analytic sheaves has a natural analytic structure:
one simply says that 
$\boldsymbol\varphi\in\bHom_{\cO}(\cE,\cA\oplus\cB)$ is 
analytic if $\pr_{\cA}\boldsymbol\varphi$ and $\pr_{\cB}\boldsymbol\varphi$ 
are both analytic.
If $\cC$ is a further analytic sheaf, an $\cO$--homomorphism 
$\theta\colon\cC|U\to\cA\oplus\cB|U$ will be analytic precisely when 
$\alpha=\pr_{\cA}\theta$ and $\beta=\pr_{\cB}\theta$ are.
Now $\theta$ is uniquely determined by $\alpha$ and $\beta$, and is denoted 
$\alpha\oplus\beta$.
Conversely, an $\cO$--homomorphism $\psi\colon\cA\oplus\cB|U\to\cC|U$ will be 
analytic if its compositions $\gamma,\delta$ with the embeddings 
$\iota_{\cA}\colon\cA\to\cA\oplus\cB$, $\iota_{\cB}\colon\cB\to\cA\oplus\cB$ are; in this case we write $\psi=\gamma+\delta$.

Clearly, $\cO^E\oplus\cO^F$ is analytically isomorphic to $\cO^{E\oplus F}$.

\subhead 3.4\endsubhead
If $\cA$ is an analytic sheaf, an analytic structure can be defined on any 
$\cO$--submodule $\cS\subset \cA$ by letting
$$
\bHom(\cE,\cS)=\bHom(\cE,\cA)\cap\bHom_{\cO}(\cE,\cS).
$$
It is straightforward that (i) and (ii) in 3.1 are satisfied.
Similarly, there is an analytic structure on $\cA/\cS$.
With the projection $\pi\colon\cA\to\cA/\cS$ one lets
$$
\bHom(\cE,\cA/\cS)=\pi_*\bHom(\cE,\cA)\subset\bHom_{\cO}(\cE,\cA/\cS).\tag3.1
$$
Again, (i) in 3.1 is straightforward to check. As for (ii), let 
$\boldsymbol\tau\in\bHom_{\cO}(\cO,\cA/\cS)_x$.
Over a neighborhood $U$ of $x$ there is
a section $a\in\Gamma(U,\cA)$ such that ${\boldsymbol\tau} 1=\pi a(x)$.
If an analytic homomorphism $\alpha\colon\cO|U\to\cA|U$ is defined by 
$\alpha\bold f=\bold f a$, its germ $\boldsymbol\alpha$ at $x$ will satisfy
$$
\boldsymbol\tau=\pi_*\boldsymbol\alpha\in\pi_*\bHom(\cO,\cA)=\bHom(\cO,\cA/\cS).
$$
In view of (3.1) therefore indeed $\bHom(\cO,\cA/\cS)=\bHom_{\cO}(\cO,\cA/\cS)$.

\subhead 3.5\endsubhead
In particular, if $\varphi\colon\cA\to\cB$ is an analytic homomorphism, $\Ker\varphi\subset\cA$, $\Im\varphi\subset\cB$, and $\Coker\varphi=\cB/\Im\varphi$ all have natural analytic structures.
Also, $\varphi$ factors through an isomorphism of $\cO$--modules
$$
\psi\colon\cA/\Ker\varphi\to\Im\varphi,\tag3.2
$$
which is easily seen to be analytic.
In spite of this, $\psi$ is not necessarily an isomorphism of analytic sheaves, as the following example shows.

Let $\Omega=X=\{0\}$, and consider a closed subspace $K$ of a Banach space $E$.
The epimorphism $E\to E/K=F$ induces an analytic epimorphism $\varphi\colon\cO^E\to\cO^F$ over $\Omega$.
In this case (3.2) becomes $\psi\colon\cO^E/\cO^K\to\cO^F$, with each plain sheaf $\cO^E,\cO^F,\cO^K$ endowed with its canonical analytic structure.
Now the image of the induced map
$$
\psi_*\colon\bHom(\cO^F,\cO^E/\cO^K)\to\bHom(\cO^F,\cO^F)
$$
contains $\id_{\cO^F}$ only if $K$ is complemented.
Otherwise $\psi_*$ is not surjective, and the inverse of $\psi$ is not analytic.

\subhead 3.6\endsubhead
Another way to express the same is that given an exact sequence
$$
0\to\cA\overset\beta\to\rightarrow\cB\overset\gamma\to\rightarrow\cC\to 0
$$
of analytic homomorphisms it is not necessarily the case that 
$\cA\approx\Im\beta$ and $\cB/\Im\beta\approx\cC$ as analytic sheaves.
These isomorphisms hold precisely when for all plain $\cE$ the induced 
sequence
$$
0\to\bHom(\cE,\cA) @>\beta_*>>\bHom(\cE,\cB) @>\gamma_*>>\bHom(\cE,\cC)\to 0
$$
is also exact.
For this reason, whenever one considers analytic homomorphisms and their 
diagrams, one should also consider the diagrams obtained by applying the 
functors $\bHom(\cE,\cdot)$.
This, at least partly, motivates the definitions in Section~4 to be presented.

\subhead 3.7\endsubhead
Any sheaf $\cA$ of $\cO$--modules can be endowed with an analytic structure by setting $\bHom(\cE,\cA)=\bHom_{\cO}(\cE,\cA)$.
In addition to this maximal structure there is also a minimal analytic structure defined as follows.
If $U\subset\Omega$ is open, let us say that a homomorphism $\alpha\colon\cE|U\to\cA|U$ is of finite type if there are a finite dimensional Banach space $F$, a plain homomorphism $\varphi\colon\cE|U\to\cO^F|U$, and a homomorphism $\beta\colon\cO^F|U\to\cA|U$ such that $\alpha=\beta\varphi$.
Germs of finite type homomorphisms form a sheaf $\bHom_{\min}(\cE,\cA)$, and the choice $\bHom(\cE,\cA)=\bHom_{\min}(\cE,\cA)$ defines an analytic structure.
Clearly, any analytic structure on $\cA$ satisfies
$$
\bHom_{\min}(\cE,\cA)\subset\bHom(\cE,\cA)\subset\bHom_{\cO}(\cE,\cA).
$$
Neither of the two extreme structures will be of any importance in the sequel, except when $\cA$ is locally finitely generated.
In this case it is often useful to endow it with the minimal analytic structure $\bHom_{\min}(\cE,\cA)$.

\subhead 3.8\endsubhead
Further examples of analytic sheaves come from analytic Fr\'echet sheaves over a finite dimensional $\Omega$, as defined in [Li].
If $\cA$ is such a sheaf, the sheaves $\bHom(\cO^E,\cA)$ of so called $AF$--homomorphisms endow $\cA$ with an analytic structure in the sense of this paper.

\head 4.\ Resolutions\endhead
\definition{4.1.\ Definition}A sequence $\cA\to\cB\to\cC$ of analytic homomorphisms over $\Omega$ is called completely exact if for each pseudoconvex $U\subset\Omega$ and each plain sheaf $\cE$ over $\Omega$ the induced sequence on sections
$$
\Gamma(U,\bHom(\cE,\cA))\to\Gamma(U,\bHom(\cE,\cB))\to\Gamma(U,\bHom(\cE,\cC))
$$
is exact.
In general, a sequence of analytic homomorphisms is completely exact if its subsequences of length three are all completely exact.
\enddefinition

Completely exact sequences are clearly exact, and an exact sequence 
$0\to\cA @>\beta>>\cB\to\cC$ is completely exact precisely when $\beta$ is an 
isomorphism between the analytic sheaves $\cA$ and $\beta(\cA)\subset\cB$.

\proclaim{Proposition}If $\cB\overset\gamma\to\rightarrow\cC\to 0$ is a 
completely exact sequence over a pseudoconvex $\Omega$, $\cF$ is a plain 
sheaf over $\Omega$, and $\varphi\colon\cF\to\cC$ is an analytic homomorphism, 
then there is an analytic homomorphism $\psi\colon\cF\to\cB$ such that 
$\varphi=\gamma\psi$.
\endproclaim

\demo{Proof}Choose $\psi$ in the preimage of $\varphi$ under the surjective homomorphism
$$
\gamma_*\colon\Gamma(\Omega,\bHom(\cF,\cB))\to\Gamma(\Omega,\bHom(\cF,\cC)).
$$
\enddemo

\definition{4.2.\ Definition}An infinite completely exact sequence
$$
\ldots @>\varphi_2>> \cF_2 @>\varphi_1>>\cF_1
@>\varphi_0>>\cB\to 0\tag4.1
$$
is called a complete resolution of $\cB$ if each $\cF_n$ is plain.
\enddefinition

We shall denote (4.1) $\cF_{\bullet}\to\cB\to 0$, or even $\cB_{\bullet}\to 0$, with $\cB_0=\cB$, $\cB_n=\cF_n$ for $n\geq 1$.
The induced sequence on sections is then written 
$\Gamma(U,\bHom(\cE,\cB_{\bullet}))\to 0$.

\proclaim{Proposition}In a complete resolution (4.1) each $\Ker\varphi_{n-1}=\Im\varphi_n$ has a complete resolution.
\endproclaim

\demo{Proof}Indeed, 
$\ldots @>\varphi_{n+1}>>\cF_{n+1} @>\varphi_n>> \Im\varphi_n\to 0$ is a complete resolution.
\enddemo

\subhead 4.3\endsubhead
A homomorphism $\beta_\bullet=(\beta_n)$ of analytic complexes $\cA_\bullet$ and $\cB_\bullet$ is called analytic if each $\beta_n\colon\cA_n\to\cB_n$ is analytic.

\proclaim{Proposition}Given a complete resolution $\cA_\bullet\to 0$ and a completely exact sequence $\cB_\bullet\to 0$ over a pseudoconvex $\Omega$, any analytic homomorphism $\varphi_0\colon\cA_0\to\cB_0$ can be extended to an analytic homomorphism $\cA_\bullet\to\cB_\bullet$.
\endproclaim

\demo{Proof}Let $\alpha_n\colon\cA_{n+1}\to\cA_n$ and $\beta_n\colon\cB_{n+1}\to\cB_n$ denote the differentials of the complexes, $\alpha_{-1}=0$, $\beta_{-1}=0$.
Suppose for $0\leq j\leq n$ analytic homomorphisms $\varphi_j$
have been constructed so that $\varphi_{j-1}\alpha_{j-1}=\beta_{j-1}\varphi_j$.
Then $\beta_{n-1}(\varphi_n\alpha_n)=\varphi_{n-1}\alpha_{n-1}\alpha_n=0$; since $\Gamma(\Omega,\bHom(\cA_{n+1},\cB_\bullet))\to 0$ is exact, there is a 
$$
\varphi_{n+1}\in\Gamma(\Omega,\bHom(\cA_{n+1},\cB_{n+1}))\qquad
\text{ with }\quad \varphi_n\alpha_n=\beta_n\varphi_{n+1}.
$$
Continuing in the same way we obtain the desired homomorphism $\varphi_\bullet\colon\cA_\bullet\to\cB_\bullet$.
\enddemo
 
\proclaim{4.4.\ Theorem}Let 
$\ldots @>\alpha_1>>\cA_1 @>\alpha_0>>\cA_0\to 0$ be a completely exact 
sequence over a pseudoconvex $\Omega$.
If $\cA_n$ has a complete resolution for $n\geq 1$, then so does $\cA_0$.
\endproclaim

\demo{Proof}First we construct a completely exact sequence 
$\cB_\bullet\to\cA_0\to 0$, where each $\cB_n$ has a complete resolution, and 
$\cB_1$ is plain.
Let $\ldots\to\cF'_n\to\cF_n @>\varphi_n>>\cA_n\to 0$ be a complete 
resolution of $\cA_n,\ n\geq 1$.
By 4.3 there are analytic homomorphisms $\psi_n$ such that the diagram
$$
\CD
&&0&&0&&0\\
&&@AAA @AAA @AAA\\
\ldots @>>> \cA_3 @>\alpha_2>> \cA_2 @>\alpha_1>>\cA_1 @>\alpha_0>>&\cA_0@>>> 0\\
&&@AA\varphi_3A @AA\varphi_2A @AA\varphi_1A\\
\ldots @>>> \cF_3 @>\psi_2>>\cF_2 @>\psi_1>>\cF_1
\endCD\tag4.2
$$
commutes.
We claim that the sequence
$$
\ldots @>\beta_3>> \cF_3\oplus\Ker\varphi_2 @>\beta_2>>\cF_2\oplus\Ker\varphi_1 
@>\beta_1>>\cF_1 @>\beta_0>>\cA_0\to 0\tag4.3
$$
is completely exact, where we write $\iota_n$ for the inclusion 
$\Ker\varphi_n\to\cF_n$ and, with 
notation introduced in 3.3,
$$
\beta_0=\alpha_0\varphi_1,\qquad \beta_1=\psi_1-\iota_1,\qquad
\beta_n=(\iota_n\oplus\psi_{n-1})(\psi_n-\iota_n),\quad n\geq 2.
$$
One checks that the $\beta_n$'s do map into the sheaves indicated in (4.3).

With a pseudoconvex $U\subset\Omega$ and a plain $\cE$, apply $\Gamma(U,\bHom(\cE,\cdot))$ to both diagrams (4.2) and (4.3).
We obtain diagrams of Abelian groups
$$
\CD
\ldots @>a_3>> A_3 @>a_2>> A_2 @>a_1>> A_1 @>a_0>> A_0@>>> 0\\
&&@AA f_3A  @AA f_2A  @AA f_1A\\
\ldots @>p_3>> F_3 @>p_2>> F_2 @>p_1>> F_1\endCD\tag4.4
$$
and
$$
\ldots @>b_3>> F_3\oplus\Ker f_2 @>b_2>> F_2\oplus\Ker f_1 @>b_1>> F_1 @> b_0>> A_0 @>>> 0;\tag4.5
$$
the first is commutative, its top row is exact, and each $f_n$ is surjective.
Here
$$
b_0=a_0f_1,\qquad b_1=p_1-i_1,\qquad b_n=(i_n\oplus p_{n-1}) (p_n-i_n),\quad 
n\geq 2,
$$
where $i_n=\Ker f_n\to F_n$ is the inclusion.
It is straightforward computation that (4.5) is a complex.
To prove (4.3) is completely exact we have to show (4.5) is exact.

First, $b_0$ is surjective because both $a_0$ and $f_1$ are.
Second, if $\xi\in\Ker b_0$ then
$$
f_1\xi\in\Ker a_0=\Im a_1=\Im a_1 f_2=\Im f_1 p_1;
$$
we used the fact that $f_2$ is surjective.
Let $f_1\xi=f_1 p_1\zeta$ and $\omega=p_1\zeta-\xi\in\Ker f_1$, so that $\xi=b_1(\zeta,\omega)\in\Im b_1$.
Third, if $n\geq 1$ and $(\xi,\eta)\in\Ker b_n=\Ker (p_n-i_n)$, then $\eta=p_n\xi$ and $0=f_n p_n\xi=a_n f_{n+1}\xi$.
Hence, as before
$$
f_{n+1}\xi\in\Im a_{n+1}=\Im a_{n+1} f_{n+2}=\Im f_{n+1} p_{n+1}.
$$
Choose $\zeta$ so that $f_{n+1}\xi=f_{n+1} p_{n+1}\zeta$, then $\omega=p_{n+1}\zeta-\xi\in\Ker f_{n+1}$.
We conclude $(\xi,\eta)=b_{n+1} (\zeta,\omega)\in\Im b_{n+1}$.
Thus (4.5) is exact and (4.3) is completely exact.
It follows from the Proposition in 4.2 that each $\cF_n\oplus\Ker\varphi_{n-1}$ has a complete resolution; hence indeed there is a completely exact sequence
$$
\ldots@>>> \cB_3 @>\beta_2>> \cB_2 @>\beta_1>>\cE_1 @>\varepsilon_0>> \cA_0 @>>> 0,
$$
where $\cE_1=\cF_1$ is plain and each $\cB_n$ has a complete resolution.

Consider the completely exact sequence
$$
\ldots @>>>\cB_3 @>\beta_2>> \cB_2 @>\beta_1>>\text{Ker}\ \varepsilon_0@>>> 0.
$$
By what we have proved there are completely exact sequences
$$
\align
&\ldots @>>> \cC_4 @>>> \cC_3 @>>> \cE_2 @>\varepsilon_1>>\Ker\varepsilon_0 @>>> 0,\qquad\text{ and so}\\
&\ldots @>>> \cC_4 @>>> \cC_3 @>>>\cE_2 @>\var_1>> \cE_1 @>\var_0>> \cA_0 @>>> 0,
\endalign
$$
with $\cE_2$ plain and each $\cC_n$ having a complete resolution.
Continuing in this way we obtain a complete resolution $\cE_\bullet\to\cA_0\to 0$.
\enddemo

\proclaim{4.5.\ Theorem (``Three lemma'')}Suppose $0\to \cA @>\beta>>\cB @>\gamma>> \cC\to 0$ is a completely exact sequence over a pseudoconvex $\Omega$.
If two among $\cA,\cB$, and $\cC$ have a complete resolution, then so does the third.
\endproclaim

\demo{Proof}We can assume that $\cA\subset\cB$ and $\beta$ is the inclusion map.
\itemitem{(a)}If $\cA$ and $\cB$ have a complete resolution, then 4.4 implies that $\cC$ also has one.

In both remaining cases $\cC$ is known to have a complete resolution $\ldots@>>>\cF @>\varphi>> \cC\to 0$.
By the Proposition in 4.1 there is a commutative diagram
$$
\CD
&&&&&&0\\
&&&&&&@AAA\\
0 @>>> \cA @>\beta>> \cB @>\gamma>> \cC @>>> 0\\
&&&&@AA\psi A @AA\varphi A @.\\
&&&&\cF @= \cF,
\endCD
$$
where $\psi\colon\cF\to\cB$ is an analytic homomorphism.
We let $\iota:\Ker\varphi\to\cF$ denote the inclusion, and 
claim that the sequence
$$
0 @>>> \Ker\varphi @>\psi\oplus\iota>> 
\cA\oplus\cF @>\beta-\psi>> \cB@>>> 0\tag4.7
$$
is completely exact.

With a pseudoconvex $U\subset\Omega$ and a plain $\cE$ apply $\Gamma(U,\bHom(\cE,\cdot))$ to (4.6) and (4.7), to obtain diagrams of Abelian groups
$$
\CD
0 @>>> A @>b>> B @>c>> C @>>> 0\\
&&&&@AA pA @AA fA\\
&&&&F@=F
\endCD
$$
and
$$
0@>>>\Ker f @> p\oplus i>> A\oplus F @> b-p >> B@>>> 0.\tag4.8
$$
The first is commutative, its top row is exact, and $f$ is surjective.
The latter is clearly a complex and exact at $\Ker f$; we have to check it is exact at the next two terms.
If $(\xi,\eta)\in\Ker (b-p)$ then $p\eta=b\xi=\xi$, whence $0=cp\eta=f\eta$.
Thus $\eta\in\Ker f$ and $(\xi,\eta)=(p\oplus i)\eta\in\Im p\oplus i$.
On the other hand, if $\zeta\in B$ then with some $\eta\in F$
$$
-c\zeta=f\eta=cp\eta,\quad\text{ i.e.},\quad \zeta+p\eta=\xi\in\Ker c=A.
$$
Thus $\zeta=\xi-p\eta\in\Im (b-p)$.
We conclude that (4.8) is exact and so (4.7) is completely exact.
Note that $\Ker\varphi$ has a complete resolution by 4.2.

(b)\ Now suppose $\cA$ too has a complete resolution.
Then $\cA\oplus\cF$ also has one and by part (a) of this proof, (4.7) implies 
that so does $\cB$.

(c)\ If $\cB$, rather than $\cA$, is known to have a complete resolution, then 
by part (b) (4.7) implies that $\cA\oplus\cF$ has a complete resolution.
In view of the completely exact sequence
$$
0@>>> \cF @>>> \cA\oplus\cF @>>>\cA@>>> 0,
$$
part (a) lets us conclude $\cA$ has a complete resolution.
\enddemo

\subhead 4.6\endsubhead
There is a clear parallel between properties of coherent sheaves and sheaves that have complete resolutions.
However, one should bear in mind that if $\cA$ and $\cB$ have complete resolutions and $\varphi\colon\cA\to\cB$ is an analytic homomorphism, then $\Ker\varphi\subset\cA$ will not necessarily have a complete resolution, not even locally; nor will $\Im\varphi\subset\cB$, see 5.3 below.
Because of this, the proof of the Three lemma is more complicated than the corresponding proof for coherent sheaves.

\head 5.\ Cohesive sheaves\endhead
\definition{5.1.\ Definition}An analytic sheaf $\cA$ over $\Omega$ is called 
cohesive if each $x\in\Omega$ has a neighborhood over which $\cA$ admits a 
complete resolution.
\enddefinition

An analytic sheaf $\cA$ that is locally isomorphic to plain sheaves (in other 
words, the sheaf of sections of a locally trivial holomorphic Banach bundle) 
is cohesive.
Indeed, if $\cE$ is a plain sheaf over an open $U\subset\Omega$ such that 
$\cE\approx \cA|U$, then 
$\ldots\to 0\to \cE\overset\approx\to\rightarrow \cA|U\to 0$ is a complete 
resolution.
Over a finite dimensional $\Omega$ there are many more examples of cohesive 
sheaves.
It can be shown that coherent sheaves with their minimal analytic 
structure discussed in 3.7 are cohesive.
We do not know whether Leiterer's Banach coherent analytic Fr\'echet sheaves, 
with their analytic structure defined in 3.8, are cohesive or not.---The 
theory developed in this paper nevertheless can be generalized so
that it includes the sheaves considered by Leiterer. In the definition
of plain sheaves $\cO^E$ one can restrict to Banach spaces $E$ in some
full subcategory of all Banach spaces, and base the notion of cohesive
sheaves on this restricted class of plain sheaves. As long as the subcategory
of Banach spaces we choose is closed under (an appropriate completion 
of) countable direct sums, all our results up to Section 9
carry over. If the subcategory consists of $L^1$ spaces of discrete
measure spaces, also known as $l^1$ spaces, the resulting 
cohesive sheaves over finite dimensional $\Omega$ will be the same as 
Banach coherent analytic Fr\'echet sheaves. This
follows from [Li, 1.3 Proposition, Theorems 2.2, 2.3, and Lemma 3.4].---

Further examples of cohesive sheaves can be constructed by the following 
theorem, an immediate consequence of 4.5.

\proclaim{Theorem}Suppose $0\to\cA\to\cB\to\cC\to 0$ is a completely exact sequence over $\Omega$.
If two among $\cA,\cB$, and $\cC$ are cohesive, then so is the third.
\endproclaim

\subhead 5.2\endsubhead
If $X'$ is another Banach space, $\Omega'\subset X'$ is open, and $f\colon\Omega\to\Omega'$ is biholomorphic, sheaves of $\cO_{\Omega'}$--modules can be pulled back by $f$ to sheaves of $\cO_\Omega$--modules.
The pullback of a plain sheaf will be plain, and the pullback of an analytic sheaf will carry a natural structure of an analytic sheaf; finally, the pullback of a cohesive sheaf will be cohesive.

\subhead 5.3\endsubhead
However, not all properties of coherent sheaves carry over to cohesive ones.
The following is an adaptation of an example of Leiterer, [Li, p.~94].

\proclaim{Example}Over $\bC$ there is an analytic endomorphism $\varphi\colon\cF\to\cF$ of a plain sheaf such that the sequence
$$
0\to\Ker\varphi\to\cF\to\Im\varphi\to 0
$$
is completely exact, but neither $\Ker\varphi$ nor $\Im\varphi$ is cohesive.
\endproclaim

Let $F=l^2,\ \cF=\cO_{\bC}^F$, and with $\zeta\in\bC$ consider $\psi(\zeta)\in\Hom(F,F)$,
$$
\psi(\zeta)(y_0,y_1,\ldots)=(y_1-\zeta y_0,y_2-\zeta y_1,\ldots),\qquad y=(y_n)\in F.
$$
In fact, $\psi\colon\bC\to\Hom(F,F)$ is holomorphic, and induces an analytic endomorphism $\varphi$ of $\cF$.
The kernel $\cK$ of $\varphi$ is supported on the disc $D=\{\zeta\in\bC\colon \zeta| < 1\}$, where it is generated by the function $h(\zeta)=(\zeta^n)^\infty_{n=0}$.

As to $\cJ=\Im\varphi$, over $D$ it agrees with $\cF$.
Indeed, since the restriction of $\psi(\zeta)$ to the hyperplane $\{y_0=0\}\subset F$ is at distance $|\zeta|<1$ to an isometric isomorphism, this restriction is invertible.
It follows that over $D$ $\varphi$ has an analytic right inverse $\cF\to\cF$, whence $\cJ|D=\cF|D$; also for any plain $\cE=\cO^E$ and open $U\subset D$ the sequence
$$
0\to\Gamma (U,\bHom(\cE,\cK))\to\Gamma(U,\bHom(\cE,\cF))@>\varphi_*>> \Gamma(U,\bHom (\cE,\cJ))\to 0\tag5.1
$$
is exact.

In fact, (5.1) is exact for all open $U\subset\bC$.
To prove this we shall need two auxiliary results.

\proclaim{Proposition~1}Let $W\subset\bC$ be open and $\sigma\colon W\to\Hom(E,F)$.
If for each $e\in E$ the function $\sigma(\cdot)e\colon W\to F$ is holomorphic, then $\sigma$ itself is holomorphic.
\endproclaim

\demo{Proof}This is the content of [M, Exercise 8E], whose solution rests on Cauchy's formula and the principle of uniform boundedness.
\enddemo

\proclaim{Proposition~2}Let $c\in\partial D$, $U\subset\bC$ a neighborhood of $c$, and $\sigma\colon U\backslash\oD\to\Hom(E,F)$ holomorphic.
If for each $e\in E$ the function $\sigma(\cdot) e$ analytically continues across $c$, then so does $\sigma$.
\endproclaim

\demo{Proof}Fix $\var> 0$ so that $c$ is the unique point on $\partial(U\backslash D)$ at distance $\leq\var$ to $(1+\var)c=c'$.
With $A\in (0,\infty)$ and $B\in (0,1/\var)$ consider 
$$
E_{AB}=\{e\in E\colon \|\sigma^{(k)}
(c') e\|_F\leq AB^k k!,\ k=0,1,\ldots\}.
$$
These are closed subsets of $E$, and the assumption, combined with Cauchy estimates, implies that their union is all of $E$.
By Baire's category theorem some $E_{AB}$ has an interior point, whence $0\in\text{ int } E_{AB}$ for some (perhaps other) $A,B$; and in fact
$$
\|\sigma^{(k)} (c') e\|_F\leq AB^k k! \|e\|_E,\quad e\in E
$$
(with yet another $A,B$).
But then this implies that the Taylor series of $\sigma$ about $c'$ has 
radius of convergence $\geq 1/B>\var$, and represents the analytic 
continuation of $\sigma$ across $c$.
\enddemo

Now we return to the analysis of the Example, and show that for every plain $\cE=\cO^E$ the sequence
$$
0\to\bHom (\cE,\cK)\to\bHom(\cE,\cF) @>\varphi_*>>\bHom(\cE,\cJ)\to 0\tag5.2
$$
is exact; exactness over $D$ we already know.
As always, the issue is whether $\varphi_*$ is surjective.
To verify this, consider $c\in\bC\backslash D$ and a germ $\boldsymbol\iota\in\bHom(\cE,\cJ)_c$ represented by a section $\iota$
of $\bHom(\cE,\cJ)$ over some connected neighborhood $U\subset\bC$ of $c$.
Let $\iota$ be induced by a holomorphic function $\theta\colon U\to\Hom(E,F)$, 
and write
$$
\theta(\zeta) e=(\theta_n (\zeta) e)^\infty_{n=0}\in F,
\qquad\zeta\in U,\ e\in E.
$$
For each $e\in E$ then $\theta(\cdot) e$ induces a section of $\cJ| U$.
Since $\cK|\bC\setminus D=0$, there is a unique $F$ valued holomorphic
function $g=(g_n)$ on a neighborhood of 
$U\setminus D$ such that 
$$
\theta(\zeta)e=\psi(\zeta) g(\zeta),
\qquad\text{ i.e., }\quad\theta_\nu(\zeta) e=g_{\nu+1}(\zeta)-\zeta g_\nu(\zeta)
$$
for $\nu=0,1,2,\ldots$ and $\zeta\in U\setminus D$.

One computes
$$
-\sum_{n\leq\nu < N}\zeta^{n-\nu-1}\theta_\nu(\zeta) e=
g_n(\zeta)-\zeta^{n-N} g_N(\zeta)
$$
for $N\geq n$, hence
$$
\lim_{N\to\infty} \bigl(-\sum_{n\leq\nu < N}
\zeta^{n-\nu-1}\theta_\nu(\zeta)e\bigr)^\infty_{n=0}=g(\zeta).\tag5.3
$$
For fixed $\zeta\in U\setminus\oD$, if we vary 
$e$, the $N$'th term on the left hand side of (5.3) represents a continuous 
linear 
operator $S_N:E\to F$, 
According to (5.3) the $S_N$ converge pointwise, so that by the principle of 
uniform boundedness $\sigma(\zeta)=\lim S_N$ is a continuous linear operator.
Hence we have a function $\sigma\colon U\backslash\oD\to\Hom(E,F)$,
$$
\sigma(\zeta) e=\bigl(-\sum^\infty_{\nu=n}\zeta^{n-\nu-1}\theta_\nu(\zeta) e\bigr)^\infty_{n=0}\in F,\quad e\in E,
$$
to which Proposition 1 applies, cf. (5.3).
We conclude $\sigma$ is holomorphic; also $\psi\sigma=\theta$ over $U\backslash\oD$.

Now if $c\not\in\partial D$ then $U\backslash\oD$ is a neighborhood of $c$.
If $c\in\partial D$ then by (5.3) $\sigma(\cdot) e$ continues across $c$ for each $e\in E$.
Therefore Proposition 2 implies that $\sigma$ itself continues across $c$.
In either case there are a neighborhood $V$ of $c$ and a holomorphic $\sigma\colon V\to\Hom (E,F)$ satisfying $\psi\sigma=\theta$ over $V$.
Passing to germs at $c$ we obtain $\boldsymbol\iota\in\Im\varphi_*$, and (5.2) is indeed exact.

Note that for any open $U\subset\bC$
$$
H^1 (U,\bHom(\cE,\cK))=H^1(U\cap D,\bHom(\cE,\cO))=H^1(U\cap D,\cO^{E^*})=0,
$$
see [Bi, Theorem 4], [Bu, p. 331], [Li, Theorem~2.3], or 9.1 below.
Hence the cohomology sequence of (5.2) gives that (5.1) is exact.

Thus we proved that $0\to\cK\to\cF\to\cJ\to 0$ is completely exact.
However, $\cK=\Ker\varphi$ is not cohesive.
Indeed, on a connected neighborhood $U$ of $1\in\bC$ any analytic 
homomorphism $\cE\to\cK|U\subset\cF|U$ of a plain sheaf $\cE$ must be 0 
on $U\backslash D$, hence on all of $U$.---Therefore $\cJ=\Im\varphi$ cannot 
be cohesive, either, by virtue of the Theorem in 5.1.

\subhead 5.4\endsubhead
Similar constructions lead to two more noteworthy examples.

\proclaim{Example}There is a cohesive sheaf $\cA$ over $\bC$ whose support has interior points yet it is not all of $\bC$.
\endproclaim

We take $F=l^2$, $\cF=\cO^F$ as in the previous example, but now we define a holomorphic $\psi\colon\bC\to\Hom(F,F)$ by
$$
\psi(\zeta) (y_0,y_1,\ldots)=(y_0,y_1-\zeta y_0,y_2-\zeta y_1,\ldots).
$$
The analytic endomorphism $\varphi\colon\cF\to\cF$ induced by $\psi$ has kernel 0 this time.
As in 5.3, one shows that the sequence $0\to\cF @>\varphi>> \Im\varphi\to 0$ is completely exact.
Thus $\Im\varphi\approx\cF$ and---by 9.2 below---$\cA=\cF/\Im\varphi$ are 
cohesive; and the support of this latter is $\bC\backslash D$.

\proclaim{5.5.\ Example}There is a plain sheaf $\cF$ over $\bC$ and two cohesive subsheaves $\cA$, $\cB\subset\cF$ such that $\cA\cap \cB$ is not cohesive.
\endproclaim

Again let $F=l^2$, $\cF=\cO^F$ and consider holomorphic maps $\psi_1,\psi_2\colon\bC\to\Hom (F,F)$
$$
\align
& \psi_1(\zeta)(y_0,y_1,\ldots)=(y_0,\zeta y_0,y_1,\zeta y_1,y_2,\zeta y_2,\ldots),\\
& \psi_2(\zeta)(y_0,y_1,\ldots)=(y_0,y_1,\zeta y_1,y_2,\zeta y_2,y_3,\ldots),
\endalign
$$
and the induced analytic endomorphisms $\varphi_1,\varphi_2$ of $\cF$.
The left inverses of $\psi_1(\zeta)$, $\psi_2(\zeta)$ given by
$$
y\mapsto (y_1,y_3,y_5,\ldots),\qquad\text{ resp. }\quad y\mapsto (y_1,y_2,y_4,y_6,\ldots)
$$
induce analytic left inverses of $\varphi_1,\varphi_2$, so that $\varphi_j$ are in fact analytic isomorphisms on their images.
In particular, $\cA=\Im\varphi_1$ and $\cB=\Im\varphi_2$ are cohesive.
However, $\cA\cap\cB$ coincides with $\Ker\varphi$ of the Example in 5.3, and is not cohesive.

\head 6.\ Resolutions and cohomology\endhead

\subhead 6.1\endsubhead
Now we turn to the following two questions.
First, how can cohesive sheaves be constructed in infinite dimensional spaces, 
beyond those that are locally isomorphic to plain sheaves?
Second, do higher cohomology groups of cohesive sheaves over pseudoconvex 
sets vanish?
Eventually we shall answer the second question in the affirmative in rather general Banach spaces, and will obtain a useful if quite particular answer to the first.
In this Section and in the next we shall concern ourselves with a special
case of the questions:\ how to recognize when a sequence of analytic 
homomorphisms is a complete resolution; and, once a complete resolution 
$\cE_\bullet\to \cA\to 0$ over a pseudoconvex $\Omega$ is granted, how to 
prove $H^q(\Omega,\cA)=0$ for $q\geq 1$?
Of course, one hopes to exploit $H^q(\Omega,\cE_n)=0$, $q\geq 1$, which is 
known in a large class of Banach spaces.
It turns out that the two issues are related, and are best treated in generality greater than that of Banach spaces, at least initially.

\subhead 6.2\endsubhead
We start with a simple result on the continuity of \v Cech cohomology groups, versions of which have been well understood and widely used in complex analysis for a long time.

\proclaim{Lemma}Let $\cK$ be a sheaf of Abelian groups over a topological space $\Omega$, and $q\geq 2$ an integer.
For $N=1,2,\ldots$, let $\frak U_N,\frak V_N$, and $\frak W_N$ be families of 
open subsets of $\Omega$, each finer than the previous one.
Assume the sequences $\{\frak U_N\}^\infty_{N=1}$, $\{\frak V_N\}^\infty_{N=1}$, and $\{\frak W_N\}^\infty_{N=1}$ are increasing, and the refinement homomorphisms
$$
H^q(\frak U_N,\cK)\to H^q(\frak V_N,\cK),\qquad H^{q-1}(\frak V_N,\cK)\to H^{q-1} (\frak W_N,\cK)
$$
are zero.
Then, writing $\frak U=\bigcup^\infty_1\frak U_N$ and 
$\frak W=\bigcup_1^\infty\frak W_N$, the refinement homomorphism 
$H^q(\frak U,\cK)\to H^q(\frak W,\cK)$ is also zero. 
\endproclaim

\demo{Proof} We shall write $\frak V=\bigcup_1^\infty\frak V_N$.
Let us introduce the following notation.
If $\frak S$ is a family of open subsets of $\Omega$, 
$(C^\bullet (\frak S),\delta)$ denotes the \v Cech complex of $\frak S$, with
values in $\cK$, and $Z^\bullet (\frak S)=\Ker\delta$ the complex of cocycles.
Given a finer family $\frak T$ and a refinement map $\frak T\to\frak S$, we 
denote the induced homomorphism 
$C^\bullet (\frak S)\to C^\bullet (\frak T)$ 
by $f\mapsto f|\frak T$.---Among the families in the Lemma there are various refinement maps:\ first of all, the inclusions 
$\frak U_N\to\frak U_M\to\frak U$, etc.~for 
$N\le M$; then refinement maps $\frak W_N\to\frak V_N\to\frak U_N$, that we 
choose to commute with the inclusions; and finally, $\frak W\to\frak U$, the 
union of the maps $\frak W_N\to\frak U_N$.
We fix these maps and the induced homomorphisms on the groups of cochains.

To prove the Lemma, let $f\in Z^q(\frak U)$ represent a cohomology class $[f]\in H^q (\frak U,\cK)$.
By the assumption, for each $N$ there is a $g_N\in C^{q-1} (\frak V)$ such that $f|\frak V_N=\delta g_N|\frak V_N$.
It follows that $(g_N-g_{N+1})|\frak V_N\in Z^{q-1} (\frak V_N)$, hence, 
again by the assumption, there is an $h_N\in C^{q-2}(\frak W)$ such 
that $(g_N-g_{N+1})|\frak W_N=\delta h_N|\frak W_N$.
If we define 
$$
k_N=g_N+\delta\sum^{N-1}_1 h_n\in C^{q-1}(\frak W),
$$
then $k_{N+1}|\frak W_N=k_N|\frak W_N$, and so there is a 
$k\in C^{q-1} (\frak W)$ such that $k|\frak W_N=k_N|\frak W_N$.
Since $f|\frak W=\delta k$, the image of $[f]$ in $H^q(\frak W,\cK)$ is 0; the cocycle $f$ being arbitrary, the Lemma is proved.
\enddemo
 
\subhead 6.3\endsubhead
Let $\frak P$ be a basis of open sets in a topological space $\Omega$.
(In subsequent applications, $\Omega$ will be an open subset of a Banach space and $\frak P$ the family of pseudoconvex subsets of $\Omega$.)

\definition{Definition}We shall say that $\frak P$ is exhaustive if, given any
$P\in\frak P$ and any open cover $\frak U$ of $P$, there is an increasing
sequence of $P_N\in\frak P$, each covered by finitely many elements of
$\frak U$, such that $\bigcup_1^\infty P_N=P$.
\enddefinition

\proclaim{6.4.\ Theorem}Let $\Omega$ be a paracompact topological space, 
$\frak P$ an exhaustive basis of open sets in $\Omega$, closed under finite intersections; let $\frak O\subset\frak P$ be another basis of open sets, also closed under finite intersections.
Consider an infinite sequence
$$
\ldots\to\cA_2 @>\alpha_1>> \cA_1 @>\alpha_0>> \cA_0\to 0\tag6.1
$$
of sheaves of Abelian groups over $\Omega$.
Assume

(i)\ the induced sequence on sections
$$
\ldots\to\Gamma (U,\cA_2)\to\Gamma (U,\cA_1)\to\Gamma (U,\cA_0)\to 0\tag6.2
$$
is exact for all $U\in\frak O$; and
\itemitem{(ii)}$H^q(U,\cA_r)=0$ for $q,r\geq 1$, $U\in\frak P$.

Then
\itemitem{(a)}(6.2) is exact for all $U\in\frak P$; and
\itemitem{(b)}$H^q(U,\cA_0)=0$ for $q\geq 1$, $U\in\frak P$.
\endproclaim

If (6.2) is exact then clearly so is (6.1).
Now given an exact sequence (6.1) and a fixed $U\subset\Omega$, the sole hypothesis $H^q(U,\cA_r)=0$ for all $q,r\geq 1$ already implies 
$H^q(U,\cA_0)=0$, $q\geq 1$, provided that either $\cA_r=0$ for some $r$, 
or $\Omega$ is finite dimensional.
In the absence of these finiteness conditions it is not hard to construct 
examples where $H^q(U,\cA_r)=0$ for $r\geq 1$ but not for $r=0$.
Since the resolutions we are working with in this paper are over infinite 
dimensional $\Omega$, and typically all $\cA_r\not= 0$, we are forced to make 
the stronger assumptions (i) and (ii).
The proof will still rest on a finiteness property, the one in the definition
of exhaustivity.
When it comes to applying the Theorem in complex analysis, the hardest is precisely to prove that the family of pseudoconvex sets is exhaustive.

\demo{Proof of (b)}Since for any $U\in\frak P$ the family 
$\frak P|U=\{P\cap U\colon P\in\frak P\}$ is an exhaustive basis of open 
sets in $U$, it will suffice to verify (b) for $U=\Omega$, under the 
assumption that $\Omega$ itself is in $\frak P$.
This is what we shall do.

We start by introducing $\cK_r=\Ker\alpha_{r-1}=\Im\alpha_r\subset\cA_r$ and 
$\cK_0=\cA_0$, that fit in short exact sequences
$$
0\to \cK_{r+1}\hookrightarrow\cA_{r+1}@>\alpha_r>>\cK_r\to 0,\qquad r\geq 0.\tag6.3
$$
By (i) the associated sequences
$$
0\to\Gamma (U,\cK_{r+1})\to\Gamma (U,\cA_{r+1})\to\Gamma(U,\cK_r)\to 0\tag6.4
$$
are also exact for $U\in\frak O$.

The heart of the matter is the following: Given an increasing sequence
of families $\frak U_N\subset\frak O$ such that 
$\frak U=\bigcup_1^\infty\frak U_N$ covers
$\Omega$, there is an increasing sequence of families $\frak V_N\subset\frak O$, each 
$\frak V_N$ finer than $\frak U_N$, such that $\bigcup_1^\infty \frak V_N$
also covers $\Omega$, and the refinement homomorphism
$$
H^q(\frak U_N,\cK_r)\to H^q(\frak V_N,\cK_r)\quad\text{ is zero for } q\ge 1.
\tag6.5
$$

To verify this, choose an increasing sequence of $P_N\in\frak P$, each
covered by a finite $\frak T_N\subset\frak U$, such that 
$\bigcup_1^\infty P_N=\Omega$. At the price of introducing empty sets
at the start and repeating some of the  $P_N$'s, we can assume that
$\frak T_N\subset\frak U_N\cap\frak T_{N+1}$, and has $\le N$ elements. 
Consider the covering
$$
\frak V_N=\{V\in\frak O\colon V\subset P_N\cap T\text{ with some }
T\in\frak T_N\}
$$
of $P_N$. Thus $\frak V_N$ is finer than $\frak U_N$.

(6.4) induces exact sequences
$$
\minCDarrowwidth{.2in}
\CD
0@>>> C^\bullet (\frak U_N,\cK_{r+1}) @>>> C^\bullet (\frak U_N,\cA_{r+1}) @>>> C^\bullet (\frak U_N,\cK_r) @>>> 0\\
& & @VVV @VVV @VVV\\
0 @>>> C^\bullet (\frak V_N,\cK_{r+1}) @>>> C^\bullet (\frak V_N,\cA_{r+1}) @>>> C^\bullet (\frak V_N,\cK_r) @>>> 0.
\endCD
$$
With vertical arrows determined by a refinement map $\frak V_N\to\frak U_N$, this diagram is commutative.
The bottom row induces an exact sequence 
$$
\aligned
\ldots @>>> H^q (\frak V_N,\cA_{r+1})& @>>> H^q (\frak V_N,\cK_r) @>c>>\\
& H^{q+1} (\frak V_N,\cK_{r+1}) @>>> H^{q+1} (\frak V_N,\cA_{r+1}) @>>>\ldots
\endaligned\tag6.6
$$
Assumption (ii) implies that $\frak V_N$ is a Leray covering of $P_N$ for the sheaf $\cA_{r+1}$, hence the first and last groups in (6.6) vanish if $q\geq 1$.
It follows that the connecting homomorphism $c$ is an isomorphism.

There are also connecting homomorphisms $H^q(\frak U_N,\cK_r)\to H^{q+1} (\frak U_N,\cK_{r+1})$, but they need not be isomorphisms.
The various connecting and refinement homomorphisms make up a commutative diagram 
$$
\minCDarrowwidth{.2in}
\CD
H^q (\frak U_N,\cK_r) @>>> H^{q+1} (\frak U_N, \cK_{r+1}) @>>> \ldots @>>> H^{q+N} (\frak U_N,\cK_{r+N})\\
@VVV & & & & @V{\rho}VV\\
H^q(\frak V_N,\cK_r) @>\approx>> H^{q+1} (\frak V_N,\cK_{r+1}) @>\approx>>\ldots @>\approx>> H^{q+N} (\frak V_N,\cK_{r+N}).
\endCD
$$
Since some refinement map $\frak V_N\to\frak U_N$ factors through 
$\frak T_N$, the map $\rho$ in the
diagram factors through $H^{q+N}(\frak T_N,\cK_{r+N})$. This latter 
group, however, is zero, because $\frak T_N$ has $\le N$ elements.
Therefore $\rho$ is the zero map, and
(6.5) follows.

Similarly, there is an increasing sequence $\frak W_N\subset\frak O$, each 
$\frak W_N$ finer than $\frak V_N$, such that $\frak W=\bigcup_1^\infty\frak W_N$
covers $\Omega$ and the refinement homomorphism
$$
H^q(\frak V_N,\cK_r) @>>> H^q (\frak W_N,\cK_r)\qquad\text{is zero for }
q\geq 1.\tag6.7
$$

Now consider an arbitrary covering $\frak U\subset\frak O$ of $\Omega$,
and let $\frak U_N=\frak U$ for $N=1,2,\ldots$. Construct 
$\frak V_N,\frak W_N,\frak W$ as above. In view of 6.2, (6.5) and (6.7)
imply $H^q (\frak U,\cK_r)\to H^q(\frak W,\cK_r)$ is zero for $q\geq 2$, 
and therefore so is the canonical homomorphism 
$H^q(\frak U,\cK_r)\to H^q (\Omega,\cK_r)$.
Since the inductive limit of these homomorphisms, as $\frak U$ ranges over 
coverings $\subset\frak O$, is the identity map of $H^q(\Omega,\cK_r)$, it 
follows that $H^q(\Omega,\cK_r)=0$ for $q\geq 2$.

To take care of $q=1$, consider a portion of the exact cohomology sequence 
associated with (6.3):
$$
H^1(\Omega,\cA_{r+1})@>>> H^1 (\Omega,\cK_r) @>>> H^2 (\Omega,\cK_{r+1}),
$$
in which the first and third terms vanish by (ii) and by what we have just proved.
It follows that the middle term also vanishes.
Since in all this proof $\Omega$ can be replaced by any $U\in\frak P$, in fact
$$
H^q (U,\cK_r)=0,\qquad q\geq 1,\quad r\geq 0,\quad U\in\frak P.\tag6.8
$$
In particular $(r=0)$ $H^q (U,\cA_0)=0$ for $q\geq 1$, as claimed.
\enddemo
\demo{Proof of (a)}With $U\in\frak P$ consider now the beginning of the exact cohomology sequence of (6.3)
$$
0\to\Gamma (U,\cK_{r+1})\to\Gamma(U,\cA_{r+1})\to\Gamma(U,\cK_r)\to H^1(U,\cK_{r+1}).\tag6.9
$$
In view of (6.8) the last term is 0, so that (6.9) is a short exact sequence; which is just another way of expressing (a).
\enddemo

\head 7.\ Resolutions and cohomology in Banach spaces\endhead
\subhead 7.1\endsubhead
Before we can apply the Theorem in 6.4 to complete resolutions over an open 
subset $\Omega$ of a Banach space, we must find a suitable family $\frak P$.
Assumption (ii) and the conclusion of the Theorem suggest that $\frak P$ 
should consist of all pseudoconvex subsets of $\Omega$.
However, it is not known whether in a general (or even separable) Banach space pseudoconvex subsets form an exhaustive family, and for this reason we shall restrict ourselves to spaces where exhaustivity is known.
Banach spaces with an unconditional basis (see 1.5) are such.
As in earlier works [L4-6], on which this one rests, we shall introduce a 
larger class of spaces $X$, and discuss cohesive sheaves in this class.
Let
$$
B_X(\rho)=B(\rho)=\{x\in X\colon \|x\| < \rho\},\qquad \rho>0,
$$
and consider the following

\proclaim{Hypothesis}
There is a $\rho\in(0,1)$ such that for any Banach space $E$, holomorphic function $f\colon B(1)\to E$, and $\var>0$, there is a holomorphic function $g\colon X\to E$ satisfying $\|f(x)-g(x)\|_E<\var$ if $x\in B(\rho)$.
\endproclaim

All Banach spaces with an unconditional basis satisfy the Hypothesis, see 
[L3,J2], and some without unconditional basis also satisfy it, see [P1,Me].
It may very well be that the Hypothesis holds in all separable 
Banach spaces.

\proclaim{7.2.\ Lemma}If $X$ satisfies the Hypothesis above then so does $X\oplus\bC$.
\endproclaim

\demo{Proof}Whether a space satisfies the Hypothesis depends only on its topology (although passing to an equivalent norm may affect the value of $\rho$).
We are therefore free to choose any norm on $X\oplus\bC$; we shall use $\|(x,\zeta)\|=\max(\|x\|,|\zeta|)$.
Given a Banach space $E$, consider the space $F$ of bounded holomorphic functions $B_{\Bbb C}(1/2)\to E$, endowed with the sup norm.
Any holomorphic function $f\colon B_{X\oplus\bC}(1)\to E$ gives rise to a holomorphic function $\varphi\colon B_X(1)\to F$,
$$
\varphi(x)=f(x,\cdot),\qquad x\in B_X(1),
$$
which, by assumption, can be approximated by holomorphic $\psi\colon X\to F$, uniformly on some $B_X(\rho)$.
Here $0<\rho<1/2$ is independent of $f$.
The function $h(x,\zeta)=\psi(x)(\zeta)$ is holomorphic on $X\times B_{\bC}(1/2)$, and, according to [L2, Th\'eor\`eme 1.1], can be approximated by a holomorphic $g\colon X\oplus\bC\to E$, uniformly on $B_{X\oplus\bC}(\rho)$. From this 
the Lemma follows.
\enddemo

\subhead 7.3\endsubhead
The Hypothesis in 7.1 concerns us because of the following two results, proved in [L5,6]:

\proclaim{Theorem~1 (Plurisubharmonic domination)}Suppose a Banach space $X$ has a Schauder basis and satisfies the Hypothesis in 7.1.
If $\Omega\subset X$ is pseudoconvex and $u\colon\Omega\to\bR$ is a locally bounded function, then there is a continuous plurisubharmonic function $v\colon\Omega\to\bR$ such that $v\geq u$.
\endproclaim

\proclaim{Theorem~2}If $X$ and $\Omega$ are as above, $L\to\Omega$ is a locally trivial holomorphic Banach bundle, and $q\geq 1$, then $H^q(\Omega,L)=0$.
In particular, if $\cA$ is a plain sheaf over $\Omega$, then $H^q(\Omega,\cA)=0$.
\endproclaim

According to [P6, Theorem 1.3], in Theorem 2 the Hypothesis from 7.1 can
be replaced by the assumption that in $\Omega$ plurisubharmonic domination
is possible (in the sense of Theorem 1).

\proclaim{7.4.\ Lemma}Suppose $X$ has a Schauder basis and the Hypothesis in 7.1 holds.
If $\Omega\subset X$ is open, then the family $\frak P$ of pseudoconvex subsets of $\Omega$ is exhaustive.
\endproclaim

\demo{Proof}Let $\frak U$ be an open cover of $P\in\frak P$. To produce
the required sequence $P_N$, we can
assume by Lindel\"of's theorem that $\frak U=\{U_1,U_2,\ldots\}$ is
countable. For $x\in P$ define $u(x)$ to be the smallest $n$ such
that $x\in U_n$.
By plurisubharmonic domination there is a plurisubharmonic $v\colon U\to\bR$ 
such that $u\leq v$; then $P$ is the increasing union of 
$P_N=\{x\in P\colon v(x)< N\}\subset\bigcup_1^NU_n$, and each $P_N\in\frak P$.
\enddemo

\proclaim{7.5.\ Corollary}If $X$ and $\Omega$ are as in 7.4,
and $\cB_\bullet\to 0$ is a complete resolution over $\Omega$, then for any 
plain sheaf $\cE$ over $\Omega$ and pseudoconvex $U\subset\Omega$
$$
H^q (U,\bHom(\cE,\cB_0))=0,\qquad q\geq 1.\tag7.1
$$
\endproclaim

\demo{Proof}Apply 6.4 with $\frak P=\frak O$ consisting of all pseudoconvex 
subsets of $\Omega$, and $\cA_r=\bHom(\cE,\cB_r)$.
The assumptions in 6.4 are satisfied by 7.4 and Theorem~2 above, and (7.1) 
follows.
\enddemo

\proclaim{7.6.\ Lemma}Let $X,\Omega$ be as in 7.4, and
$$
\cB_\bullet\to\cB_0\to 0\tag7.2
$$
a sequence of analytic homomorphisms over $\Omega$, with $\cB_r$ plain for 
$r\geq 1$.
If each $x\in\Omega$ has a neighborhood over which (7.2) is completely 
exact, then it is completely exact over $\Omega$.
\endproclaim

\demo{Proof}This follows from 6.4 if we let $\frak P$ consist of all 
pseudoconvex subsets of $\Omega$, $\frak O$ consist of those $U\in\frak P$ 
over which (7.2) is completely exact, and $\cA_r=\bHom(\cE,\cB_r)$, $\cE$ an 
arbitrary plain sheaf.
Again, the hypotheses in 6.4 are satisfied by Theorem~2 and 7.4 above.
\enddemo

\head 8.\ Gluing Complete Resolutions\endhead
\subhead 8.1\endsubhead
In this Section $X$ will denote a Banach space with a Schauder basis that 
satisfies the Hypothesis in 7.1, and $\Omega\subset X$ will be pseudoconvex.
We shall prove a few results to the effect that if an analytic sheaf $\cA$ 
over $\Omega$ has complete resolutions over certain subsets of $\Omega$, 
then it has one over all of $\Omega$.
The corresponding result in finite dimensions depends on a lemma of Cartan 
concerning holomorphic matrices, which is a way of saying that certain 
holomorphic vector bundles are trivial.
We start with a similar result in our infinite dimensional setting.

\proclaim{8.2.\ Lemma}Let $\pi\colon L\to\Omega\times\bC$ be a locally 
trivial holomorphic Banach bundle.
If the restriction $L|\Omega\times\{\zeta\}$ is trivial for some 
$\zeta\in\bC$, then it is trivial for all $\zeta\in\bC$.
\endproclaim

This also follows easily from [P6, Theorem 1.3d].

\demo{Proof}First we show that $L$ admits a holomorphic connection, i.e.,
there is a holomorphic subbundle of $TL$ that is complementary to 
$\Ker\pi_*$, and is compatible with the vector bundle operations.
If over an open $U\subset\Omega\times\bC$ there is a holomorphic connection $H\subset TL|U$, and $\theta$ is an $L$--valued holomorphic 1--form on $U$, another holomorphic connection $H^\theta$ over $U$ can be constructed as 
follows. For
$y\in U$, $l\in L_y$, and $\eta\in T_y U$, let $\theta_l (\eta)\in T_l L_y$ 
denote the vector corresponding to $\theta(\eta)\in L_y$ under the canonical isomorphism $L_y\approx T_l L_y$.
Then define
$$
H_l^\theta=\{\lambda\in T_l L\, \colon\, \lambda-\theta_l (\pi_*\lambda)\in H_l\},\qquad\text{ and }\quad H^\theta=\bigcup_{l\in L|U} H_l^\theta.
$$
One easily checks that the subbundle $H^\theta\subset TL|U$ thus obtained is indeed a connection; conversely, any holomorphic connection $H'$ on $L|U$ is obtained from $H$ via a unique holomorphic 1--form $\theta$; and $H^{\theta_1+\theta_2}=(H^{\theta_1})^{\theta_2}$.

Since $L$ is locally trivial, $\Omega\times\bC$ can be covered by open subsets $U$ over which $L$ admits a holomorphic connection $H_U$.
If $U\cap V\neq\emptyset$, there is an $L$--valued holomorphic 1--form 
$\theta_{UV}=\theta$ on $U\cap V$ such that $H_V=H_U^\theta$ over $U\cap V$.
The $\theta_{UV}$ form a holomorphic cocycle with values in the bundle of $L$--valued 1--forms.
Now $X\oplus\bC$ has a Schauder basis and, by 7.2, satisfies the Hypothesis in 7.1.
Therefore Theorem~2 in 7.3 applies, and it follows that there are holomorphic $L$--valued 1--forms $\theta_U$ on each $U$ such that $\theta_U-\theta_V=\theta_{UV}$.
This latter means that over $U\cap V$ the connections $H_U^{\theta_U}$ and 
$H_V^{\theta_V}$ coincide, and so define a holomorphic connection $H$ on $L$.

Given $H$, we can construct an isomorphism between any two restrictions $L|\Omega\times\{\zeta\},\ L|\Omega\times\{\zeta'\}$.
If $l\in L_{(x,\zeta)}$, lift the line segment joining $(x,\zeta)$ with $(x,\zeta')$ to a curve in $L$, starting at $l$ and everywhere tangential to $H$.
Denoting the endpoint of the curve $\Phi(l)$, the map $\Phi$ will be an isomorphism between the two restricted bundles.
Therefore if one of them is trivial, so is the other.
\enddemo

\subhead 8.3\endsubhead
For purposes of this Section we introduce the following terminology.
Let $U\subset X$ be open, $E$ a Banach space, and GL$(E)\subset\Hom (E,E)$ the open subset of invertible homomorphisms.
This is in fact a Banach--Lie group.

\definition{Definition}A holomorphic map $\varphi\colon U\to\text{ GL}(E)$ is called connectible if there is a holomorphic $\tilde\varphi\colon U\times\bC\to \text{ GL}(E)$ such that $\tilde\varphi(x,0)=\id_E$ and $\tilde\varphi(x,1)=\varphi(x)$, $x\in U$.
The map $\tilde\varphi$ is called a connecting map of $\varphi$.
\enddefinition

\proclaim{Proposition}Let $\cA$ be an analytic sheaf over a pseudoconvex $W\subset X$.
Given Banach spaces $E,E'$, and completely exact sequences
$$
\cO^E|W @>\alpha>> \cA @>>> 0\qquad\text{ and }\qquad
\cO^{E'}|W@>\alpha'>>\cA\to 0,
$$
denote by $\pi,\pi'$ the projections $\cO^{E\oplus E'}\to\cO^E$, resp.~$\cO^{E'}$; then there is a connectible $\var\colon W\to\text{ GL}(E\oplus E')$ such that $\alpha'\pi'=\alpha\pi\var$.
\endproclaim

Above we identified $\var$ with the analytic endomorphism of $\cO^{E\oplus E'}_W$ it induces.

\demo{Proof}By the Proposition in 4.1 the homomorphisms $\alpha'$ and $\alpha$ can be factored
$$
\alpha'=\alpha\psi,\qquad \alpha=\alpha'\varphi,
$$
with analytic homomorphisms
$$
\psi\colon \cO^{E'}|W\to\cO^E|W,\qquad \varphi\colon \cO^E|W\to\cO^{E'}|W
$$
Again, we shall denote the corresponding holomorphic maps $W\to\Hom (E',E)$, resp.~$\Hom(E,E')$ by the same symbols $\psi,\varphi$.
Defining $\tilde\var\colon W\times\bC\to\text{ GL}(E\oplus E')$ by
$$
\tilde\var (x,\zeta) \pmatrix e\\ e'\endpmatrix=\pmatrix \id_E&\zeta\psi(x)\\ 0&\id_{E'}\endpmatrix
\pmatrix \id_E&0\\ \zeta\varphi(x)&\id_{E'}\endpmatrix^{-1}\pmatrix e\\ e'\endpmatrix,\qquad\pmatrix e\\ e'\endpmatrix\in E\oplus E',
$$
it is straightforward that $\var(x)=\tilde\var (x,1)$ does it.
\enddemo

\proclaim{8.4.\ Proposition}Let $V,V'\subset X$ be open, $V\cap V'=W$ and $V\cup V'$ pseudoconvex.
If an analytic sheaf $\cA$ over $V\cup V'$ has complete resolutions over $V$ and $V'$ then it has one over $V\cup V'$.
\endproclaim

\demo{Proof}Consider the end portion of complete resolutions over $V$ and $V'$
$$
\cO^E|V @>\alpha>> \cA|V\to 0\qquad\text{ and }\qquad\cO^{E'}|V' @>\alpha'>> \cA|V'\to 0.\tag8.1
$$
Apply the previous Proposition with the restrictions $\alpha|W$, $\alpha'|W$.
Putting $F=E\oplus E'$, there are a holomorphic $\var\colon W\to\text{ GL}(F)$ 
and a connecting map $\tilde\var\colon W\times\bC\to\text{ GL}(F)$ such that
$$
\alpha'\pi'=\alpha\pi\var\qquad\text{ over }W.\tag8.2
$$
Glue the trivial bundles
$$
F\times (V\times\bC)\to V\times\bC\qquad\text{ and }\qquad F\times (V'\times\bC)\to V'\times\bC
$$
together with the gluing map $\tilde\var$, to obtain a holomorphic Banach bundle $L\to (V\cup V')\times\bC$.
Since $\tilde\var(x,0)=\id,\ L|(V\cup V')\times \{0\}$ is trivial; by 8.2 $L$ is therefore trivial over $(V\cup V')\times \{1\}$ as well.
Hence there are holomorphic
$$
\delta\colon V\to\text{ GL}(F)\qquad\text{ and }\qquad\delta'\colon V'\to\text{ GL}(F)
$$
such that $\var=\delta\delta^{'-1}$.
Since the analytic homomorphisms
$$
\alpha\pi\delta\colon \cO^F|V\to\cA|V\qquad\text{ and }
\qquad\alpha'\pi'\delta'\colon\cO^F|V'\to \cA|V'
$$
agree over $V\cap V'$ according to (8.2), they induce an analytic homomorphism 
$\varphi_0\colon\cO^F|V\cup V'\to\cA$.

Call $F=F_1$, and let $\cA_1=\Ker\varphi_0$.
The sequence $0\to\cA_1\hookrightarrow\cO^{F_1}@>\varphi_0>>\cA\to 0$ is 
completely exact over $V$ and $V'$, because the sequences (8.1) were; and by 
the Three lemma in 4.5 both $\cA_1|V$ and $\cA_1|V'$ have complete resolutions.
We can therefore repeat our construction above, and obtain sequences of analytic homomorphisms
$$
0\to\cA_r\hookrightarrow \cO^{F_r} @>\varphi_{r-1}>> \cA_{r-1}\to 0
$$
over $V\cup V'$, that are completely exact over $V$ and $V'$.
These short exact sequences can be consolidated in a sequence
$$
\ldots\to \cO^{F_2} @>\varphi_1>> \cO^{F_1} @>\varphi_0>> \cA\to 0,
$$
itself completely exact over $V$ and $V'$.
In view of 7.6 it is therefore a complete resolution over $V\cup V'$.
\enddemo

\proclaim{8.5.\ Proposition}Suppose $\Omega$ is the increasing union of pseudoconvex $\Omega_N\subset\Omega$, $N\in\bN$, and $\cA$ is an analytic sheaf over $\Omega$.
If $\cA$ has a complete resolution over each $\Omega_N$ then it has one over $\Omega$ too.
\endproclaim

\demo{Proof}There are completely exact sequences
$$
\cO^{E_N}|\Omega_N @>\alpha_N>> \cA|\Omega_N\to 0,
$$
with Banach spaces $(E_N, \|\ \|_N)$.
Consider the $l^\infty$ sum
$$
F=\{ e=(e_N)^\infty_1\colon\text{ each }e_N\in E_N\text{ and }\|e\|_F=\sup_N \|e_N\|_N < \infty\},
$$
and the projections $\pi_N\colon\cO^F\to\cO^{E_N}$.

For each $N$ apply the Proposition in 8.3 with $\alpha=\alpha_N$, $\alpha'=\alpha_{N+1}|\Omega_N$.
The resulting connectible $\var\colon\Omega_N\to\text{GL}(E_N\oplus E_{N+1})$
can be extended to a holomorphic $\var_N\colon\Omega_N\to\text{GL}(F)$, by 
letting $\var_N(x)$ for $x\in\Omega_N$ act by identity on each $E_n$, $n\neq N, N+1$.
Thus we have
$$
\alpha_{N+1}\pi_{N+1}=\alpha_N\pi_N\var_N\qquad\text{over }\Omega_N.\tag8.3
$$
The maps $\var_N$ are also connectible, let $\tilde\var_N$ be their connecting map.

Next construct a holomorphic Banach bundle $L\to\Omega\times\bC$ with fiber $F$, by gluing together trivial bundles over $\Omega_N\times\bC$ with gluing maps
$$
\tilde\var_{M-1}\tilde\var_{M-2}\ldots\tilde\var_N\colon (\Omega_M\cap\Omega_N)\times\bC\to\text{ GL}(F),
$$
$M>N$.
Again $L|\Omega\times\{0\}$ is trivial, whence so is $L|\Omega\times\{1\}$ by 
8.2.
Therefore there are holomorphic $\delta_N\colon\Omega_N\to\text{ GL}(F)$ such 
that $\var_N=\delta_N\delta^{-1}_{N+1}$.
Now (8.3) implies the maps $\alpha_{N+1}\pi_{N+1}\delta_{N+1}$ and 
$\alpha_N\pi_N\delta_N$ agree over $\Omega_N$, hence define an analytic
homomorphism $\varphi_0\colon\cO^F\to\cA$.

Write $F_1$ for $F$ and $\cA_1=\Ker\varphi_0$, to obtain, as in 8.4, a 
sequence $0\to\cA_1\hookrightarrow\cO^{F_1}\to\cA\to 0$, completely exact over 
each $\Omega_N$.
By the Three lemma in 4.5, each $\cA_1|\Omega_N$ has a complete resolution.
As in 8.4, we can repeat our construction to produce sequences of analytic homomorphisms
$$
0\to \cA_r\hookrightarrow \cO^{F_r}@>\varphi_{r-1}>> \cA_{r-1}\to 0,
$$
that are completely exact over each $\Omega_N$, and give rise to a sequence
$$
\ldots\to\cO^{F_2}@>\varphi_1>>\cO^{F_1}@>\varphi_0>>\cA\to 0,\tag8.4
$$
itself completely exact over each $\Omega_N$.
By 7.6 we can conclude that (8.4) is in fact a complete resolution.
\enddemo

\subhead 8.6\endsubhead
The Proposition in 8.4 is sufficient to construct a complete resolution of a 
cohesive sheaf over compact subsets of a pseudoconvex $\Omega$, and this is 
all Cartan needed in the finite dimensional case he studied.
The more precise result of Leiterer, that a (complete) resolution exists over 
all of $\Omega$---still assuming $\dim\Omega < \infty$---depends on 8.5 as well.
To deal with infinite dimensional $\Omega$ one more ingredient will be needed.

\proclaim{Proposition}Let $\cA$ be an analytic sheaf over $\Omega$; let $n\in\bN$, $\tilde p\in\Hom(X,\bC^n)$, and $p=\tilde p|\Omega$.
If $D=p(\Omega)\subset\bC^n$ is pseudoconvex and can be covered by open $V\subset D$ such that $\cA|p^{-1} V$ has a complete resolution, then $\cA$ itself has a complete resolution.
\endproclaim

\demo{Proof}We first show that $\cA|p^{-1} U$ has a complete resolution if $U\subset D$ is relatively compact and pseudoconvex.
Suppose for some $U$ it does not.
Let $s\colon\bC^n\to\bR$ be $\bR$--linear, $c<d$ real numbers, and
$$
U^+=\{z\in U\colon s(z) > c\},\qquad U^-=\{z\in U\colon s(z) < d\}.
$$
Both $U^\pm$ are pseudoconvex and $U^+\cup U^-=U$.
It follows from 8.4 that $\cA$ cannot have a complete resolution over both
$p^{-1} U^\pm$.
Denote by $U_1$ one of $U^\pm$ so that $\cA$ has no complete resolution over 
$p^{-1} U_1$.
Choosing further linear forms $s_1,s_2,\ldots$ and $c_1<d_1,c_2<d_2,\ldots$, 
etc.~we obtain a sequence $U\supset U_1\supset U_2\supset\ldots$ with the 
property that $\cA|p^{-1} U_k$ has no complete resolution, $k\in\bN$.
A judicious choice of $s_k,c_k,d_k$ will ensure that $\diam\ U_k\to 0$, so 
that the $U_k$ converge to some $z\in D$.
Hence the assumption implies that $\cA|p^{-1} U_k$ has a complete resolution 
for some $k$, after all; this contradiction shows that $\cA|p^{-1} U$ 
must have 
a complete resolution.

Once this granted, exhaust $D$ by an increasing sequence of relatively 
compact pseudoconvex $D_N\subset D$.
Since $\Omega$ is the increasing union of the pseudoconvex subsets 
$p^{-1} D_N$, the Proposition follows from 8.5.
\enddemo

\head 9.\ The Main Theorem\endhead
\proclaim{9.1.\ Theorem}Suppose a Banach space $X$ has a Schauder basis and satisfies the Hypothesis in 7.1.
If $\cA$ is a cohesive sheaf over a pseudoconvex $\Omega\subset X$, then
\itemitem{(a)}$\cA$ admits a complete resolution;
\itemitem{(b)}$H^q(\Omega,\cA)=0$ for $q\geq 1$; and more generally,
\itemitem{(c)}$H^q(\Omega,\bHom(\cE,\cA))=0$ for $q\geq 1$ and $\cE$ a plain sheaf.
\endproclaim

Relying on [P6, Theorem 1.3], here the Hypothesis from 7.1 can
be replaced by assuming that in every pseudoconvex subset of $\Omega$
plurisubharmonic domination is possible; see the remark in 7.3.

\demo{Proof}Fix a Schauder basis $e_1,e_2,\ldots\in X$, and let $\pi_N\colon X\to X$ be the projection
$$
\pi_N(\sum^\infty_1\lambda_n e_n)=\sum^N_1\lambda_n e_n.
$$
As explained in [L4, Section~7], it can be arranged that the norm $\|\ \|$ in 
$X$ is such that all projections $\pi_N-\pi_M$ have norm $\leq 1$, and in what 
follows we shall assume this.
We set $B(x,R)=\{y\in X\colon \|x-y\| < R\}$.

There is a function $u\colon \Omega\to\bR$ such that for $x\in\Omega$ the restriction of $\cA$ to $B(x,2 e^{-u(x)})\subset\Omega$ has a complete resolution; in fact, $u$ can be chosen locally bounded.
Plurisubharmonic domination in 7.3 yields a continuous plurisubharmonic function $v\colon\Omega\to (0,\infty)$ such that $v\geq u$.
Thus $\cA|B(x,2e^{-v(x)})$ has a complete resolution when $x\in\Omega$.
Define
$$
\gather
D_N=\{z\in\Omega\cap\pi_N X\,\colon\, e^{v(z)} < 2N+2\},\\
\Omega_N=\{x\in\pi_N^{-1} D_N\,\colon\, \| x-\pi_N x\| < e^{-v(\pi_N x)}\}.
\endgather
$$

According to [L6, Proposition 4.3], each $\Omega_N\subset\Omega$ is pseudoconvex.
Clearly, each $z\in D_N$ has a neighborhood $V\subset D_N$ such that
$$
\Omega_N\cap \pi_N^{-1} V\subset B (z, 2e^{-v(z)});
$$
in particular, $\cA|\Omega_N\cap\pi_N^{-1} V$ has a complete resolution.
By the Proposition in 8.6 therefore $\cA|\Omega_N$ also has a complete resolution.

To conclude, note that according to [L6, Proposition 4.3] $\Omega'_N=\bigcap_{n\geq N}\Omega_n$ is a locally finite intersection---hence each $\Omega'_N$ is (open and) pseudoconvex---and $\Omega=\bigcup\Omega'_N$.
Since each $\cA|\Omega'_N$ also has a complete resolution, 8.5 applies and we obtain (a).
Finally, (b) and (c) follow from (a) by virtue of 7.5.
\enddemo

\proclaim{9.2.\ Lemma}Suppose a Banach space $X$ has a Schauder basis and 
satisfies the Hypothesis in 7.1. Let
$\cA$ be an analytic sheaf over an open $\Omega\subset X$ and 
$\cB\subset\cA$ a cohesive subsheaf. If one of  $\cA$ and $\cA/\cB$ is 
cohesive, then so is the other.
\endproclaim

\demo{Proof}The way subsheaves and quotient sheaves were endowed with an 
analytic structure in 3.4 implies that for any plain sheaf $\cE$ over 
$\Omega$
$$
0\to\bHom(\cE,\cB)\to\bHom(\cE,\cA)\to\bHom(\cE,\cA/\cB)\to 0
$$
is exact.
If $U\subset\Omega$ is pseudoconvex, the associated cohomological sequence gives
$$
\aligned
0\to\Gamma(U,\bHom(\cE,\cB))&\to\Gamma(U,\bHom(\cE,\cA))\to\\
&\Gamma (U,\bHom(\cE,\cA/\cB))\to H^1 (U,\bHom(\cE,\cB))=0
\endaligned
$$
in view of 9.1.
Hence $0\to \cB\to\cA\to \cA/\cB\to 0$ is completely exact and 
the claim follows from the Theorem in 5.1.
\enddemo

\head 10.\ Sheaves Associated with Submanifolds\endhead
\subhead 10.1\endsubhead
Consider Banach spaces $X$ and $F$, a complemented subspace $Y\subset X$, 
$k\in\bN$, and the sheaf $\cJ\subset\cO^F$ of germs vanishing on $Y$ to order 
$k$.
Thus over $X\setminus Y$ the sheaves $\cJ$ and $\cO^F$ coincide.
There is a complex of analytic homomorphisms
$$
\ldots\to\cH_2 @>\varphi_1>> \cH_1 @>\varphi_0>> \cH_0=\cJ\to 0,\tag10.1
$$
that generalizes Koszul's complex (to which it reduces when codim $Y<\infty$ and $k=1$), defined as follows.
Fix a closed complement $Z\subset X$ of $Y$, and let $E_r$ denote the Banach space of alternating $r$--linear forms $Z\oplus Z\oplus\ldots\oplus Z\to F$, $E_0=F$.
In the natural identification of the tangent spaces $T_x X$ with $X$, the subspaces $T_x(Z+x)\subset T_x X$ correspond to $Z$.
Hence sections of $\cO^{E_r}$ over $U\subset X$ can be identified with holomorphic {\sl relative} $r$--forms $e$ on $U$, with values in $F$; relative meaning that $e$ can be evaluated only on $r$--tuples $\xi_j\in T_x X,j=1,\ldots,r$, that are tangent to $Z+x$.
We let $\cH_r\subset\cO^{E_r}$ denote the sheaf of such holomorphic relative $r$--forms, which in addition vanish on $Y$ to order $k-r$ if $r<k$.
Thus $\cH_r=\cO^{E_r}$ when $r\geq k$, and $\cH_0=\cJ$.

Furthermore, consider the flow $g_t$ on $X$ given by
$$
g_t(y+z)=y+(\exp t)z,\qquad y\in Y,\ z\in Z,\ t\in\bR.
$$
Its infinitesimal generator $\xi=dg_t/dt|_{t=0}$ is a vector field on $X$, 
tangential to the subspaces $Z+x$.
In the identification of $T_{y+z} X$ with $X$, the vector $\xi(y+z)$ 
corresponds to $z\in Z\subset X$.
Contraction with $\xi$ defines an analytic homomorphism 
$\cO^{E_{r+1}}\ni\bold e\mapsto \iota_\xi\bold e\in\cO^{E_r}$.
We set $\varphi_r=\iota_\xi|\cH_{r+1}$ to construct the complex (10.1).

\proclaim{10.2.\ Lemma}Let $V\subset Y$, $W\subset Z$ be open and convex.
If $U=V\times W\subset X$ then the induced complex 
$\Gamma(U,\cH_\bullet)\to 0$ is exact.
\endproclaim

A version of the Lemma in the case $k=1$ was also found, with the same
proof, by Patyi, see [P7, Proposition 2.1].

\demo{Proof}Since the flow $g_t$ preserves each submanifold $Z+x$, it induces 
a pullback operation $g_t^*$ on the sheaf $\cO^{E_r}$ of relative $r$--forms, 
and one can define a Lie derivative
$$
L_\xi=dg_t^*/dt|_{t=0}\colon\cO^{E_r}\to\cO^{E_r}.
$$
There are also relative exterior derivatives $d\colon\cO^{E_r}\to\cO^{E_{r+1}}$.
Both $L_\xi$ and $d$ are homomorphisms of sheaves of Abelian groups and map the complex $\cH_\bullet\to 0$ into itself.
Note that E.~Cartan's identity 
$$
L_\xi=\iota_\xi d+d\iota_\xi\tag10.2
$$
holds (since it holds on finite dimensional $g_t$--invariant subspaces of each $Z+x$).

To prove the Lemma we have to distinguish between two cases.
If $0\not\in W$, by the Banach--Hahn theorem there is a $p\in\Hom(Z,\bC)$ such that $p(z)\neq 0$ for $z\in W$.
We extend $p$ to $X$, constant on subspaces $Y+x$.
Using the identification $T_x(Z+x)\approx Z$, $p$ induces a 1--form on each $T_x(Z+x)$, hence a holomorphic relative 1--form $\omega$ on $X$.
Clearly $\iota_\xi\omega=p$.
Suppose $h\in\Gamma (U,\cH_{r+1})$ is in the kernel of $\varphi_{r*}$; thus $\iota_\xi h=0$.
Setting $f=\omega\wedge h/p\in\Gamma(U,\cH_{r+2})$, we compute
$$
i_\xi f=(i_\xi\omega)\wedge h/p=h,
$$
so that $h\in\Im\varphi_{r+1*}$, as needed.

On the other hand, if $0\in W$ then $g_t U\subset U$ for $t\leq 0$.
Using (10.2) we write for $h\in\Ker\varphi_{r*}$
$$
h=\int^0_{-\infty} {d\over dt} (g^*_t h) dt=\int^0_{-\infty} (L_\xi g_t^* h) dt=i_\xi\int^0_{-\infty} (g_t^* dh) dt,
$$
with the improper integrals converging exponentially fast.
It follows that
$$
f=\int^0_{-\infty} (g_t^* dh) dt\in\Gamma(U,\cH_{r+2}),
$$
and $h=\varphi_{r+1*}f\in\Im \varphi_{r+1*}$, which proves the claim.
\enddemo

\proclaim{10.3.\ Theorem}Suppose $X$ is a Banach space that has a Schauder basis and satisfies the Hypothesis in 7.1.
If $F$ is another Banach space, $\Omega\subset X$ is open, $M\subset\Omega$ a direct submanifold (see~1.3), and $k\in\bN$, then the sheaf $\cJ\subset\cO^F$ of germs vanishing on $M$ to order $k$ is cohesive.
\endproclaim

\demo{Proof}We prove the Theorem by induction on $k$, so we assume it is true when $k$ is replaced by any $j$, $1\leq j < k$.
Since cohesion is a local property, and invariant under biholomorphisms, we can take $\Omega=X$ and $M=Y$ a complemented subspace with complement $Z$.
We shall show that the sequence (10.1) constructed above is completely exact.
To this end we apply 6.4, with $\frak P$ the collection of pseudoconvex subsets of $\Omega=X$, $\frak O$ the collection of $V\times W$, where $V\subset Y$, $W\subset Z$ are open and convex, and $\cA_\bullet=\cH_\bullet$.
According to 7.4 $\frak P$ is exhaustive, and assumption (i) of 6.4 has just been verified.
Since by the inductive hypothesis $\cH_r$ is cohesive for $r\geq 1$, assumption (ii) is also satisfied in view of 9.1.
Therefore 6.4 gives that $\Gamma(U,\cH_\bullet)\to 0$ is exact for $U\in\frak P$.

More generally, if $E$ is an arbitrary Banach space, the complex $\bHom(\cO^E,\cH_\bullet)\to 0$ is isomorphic to a complex of type (10.1), but constructed with $\Hom(E,F)$ replacing $F$.
Hence by what has already been proved, $\Gamma(U,\bHom(\cO^E,\cH_\bullet))\to 0$ is exact for all $U\in\frak P$, i.e.~(10.1) itself is completely exact.

Now we conclude by 9.1 and 4.4.
By the former, and by the inductive hypothesis, $\cH_r$ for $r\geq 1$ have complete resolutions.
By the latter, this implies $\cH_0=\cJ$ also has a complete resolution, in particular it is cohesive.
\enddemo

\proclaim{10.4.\ Corollary}Let $X,F,\Omega$, and $M$ be as in 10.3.
If, in addition, $\Omega$ is pseudoconvex, then any holomorphic function 
$M\to F$ is the restriction of a holomorphic function $\Omega\to F$.
\endproclaim

\demo{Proof}Let $\cJ\subset\cO^F$ denote the sheaf of germs vanishing on $M$ 
and $\cS$ the sheaf over $\Omega$ associated with the presheaf

$$
U\mapsto\{f\colon M\cap U\to F\text{ is holomorphic}\},\quad\text
{ resp. }\{0\},
$$
depending on whether $M\cap U\neq\emptyset$ or $=\emptyset$.
Restriction to $M$ defines an epimorphism of $\cO$--modules $\cO^F\to\cS$ that fits in an exact sequence $0\to\cJ\hookrightarrow\cO^F @>r>>\cS\to 0$.
By the associated exact sequence in cohomology
$$
\ldots\to\Gamma(\Omega,\cO^F) @>r_*>> \Gamma(\Omega,\cS)\to H^1(\Omega,\cJ)=0,
$$
cf.~9.1 and 10.3, $r_*$ is surjective, which is equivalent to the claim.
\enddemo

\subhead 10.5\endsubhead
It has been known for quite a while that the above Corollary fails in some Banach spaces.
Indeed, according Dineen, see [Di] and also [J1], in the nonseparable space $X=l^\infty$ there is a discrete sequence $S$ on which all holomorphic functions $X\to\bC$ are bounded.
Obviously, then, $S$ is a direct submanifold of $X$, any function $S\to\bC$ is holomorphic, but no unbounded $S\to\bC$ can be continued to a holomorphic function on $X$.

Thus, both 9.1 and 10.3 cannot generalize to $l^\infty$.
It would be of some interest to know which of the two fails in $l^\infty$.
Is it that cohesive sheaves---perhaps even plain sheaves---may have nonzero higher cohomology groups; or only that ideal sheaves of direct submanifolds are not necessarily cohesive?

However, it would be of even greater importance to clarify for what subsets $M\subset\Omega$ other than direct submanifolds will 10.3 generalize (if all other assumptions are kept).
In the Appendix we show it does not generalize to arbitrary analytic subsets.

\head 11.\ Subvarieties\endhead
\subhead 11.1\endsubhead
Cohesive sheaves can be introduced on any topological space, once we specify 
an (enriched) category of sheaves and homomorphisms to play the role of plain 
sheaves and homomorphisms.
Whether cohesive sheaves have useful properties of course depends on what 
properties the ``plain'' category specified has.
In this Section we shall study this generalization of cohesion in a context 
that is closely related to the one considered heretofore:\ for sheaves over 
subvarieties in Banach spaces.
To define subvarieties, consider a Banach space $X$ and an open 
$\Omega\subset X$.

\definition{Definition}An ideal structure on $\Omega$ is the choice of an $\cO$--submodule $\cJ^E\subset\cO^E$, for each Banach space $E$, such that for every $x\in\Omega$ and $\boldsymbol\varphi\in\cO_x^{\Hom(E,F)}$ the induced homomorphism $\boldsymbol\varphi_*$ maps $\cJ_x^E$ in $\cJ_x^F$.
The ideal structure is cohesive if each $\cJ^E$ is cohesive.
\enddefinition

\subhead 11.2\endsubhead
If $E,F,G$ are Banach spaces, there is a homomorphism
$$
\cO^{\Hom(E,F)}\otimes_{\cO} \cO^{\Hom(F,G)}\to \cO^{\Hom(E,G)}\tag11.1
$$
sending $\boldsymbol\varphi\otimes\boldsymbol\psi$ to $\boldsymbol\psi\boldsymbol\varphi$.
For any ideal structure we have

\proclaim{Proposition}(a)\ If either $\boldsymbol\varphi\in\cJ^{\Hom(E,F)}$ or 
$\boldsymbol\psi\in\cJ^{\Hom(F,G)}$ then (11.1) sends 
$\boldsymbol\varphi\otimes\boldsymbol\psi$ in $\cJ^{\Hom(E,G)}$.

(b)\ Under the obvious isomorphism $\cO^E\oplus\cO^F\to\cO^{E\oplus F}$ the 
image of $\cJ^E\oplus\cJ^F$ is $\cJ^{F\oplus G}$.
\endproclaim

\demo{Proof}(a)\ Any $\boldsymbol\varphi\in\cO_x^{\Hom(E,F)}$ induces, by 
composition, a germ $\tilde{\boldsymbol\varphi}$ of a plain homomorphism 
$\cO^{\Hom(F,G)}\to\cO^{\Hom(E,G)}$ at $x$.
According to the definition therefore
$$
\boldsymbol\psi\boldsymbol\varphi=
\tilde{\boldsymbol\varphi}\boldsymbol\psi\in\cJ_x^{\Hom(E,G)},\quad\text{if }
\boldsymbol\psi\in\cJ_x^{\Hom(F,G)}.
$$

The case when $\boldsymbol\varphi\in\cJ^{\Hom(E,F)}$ and 
$\boldsymbol\psi\in\cO^{\Hom(F,G)}$ is proved similarly.

(b)\ Applying the condition in the definition with constant germs 
$\boldsymbol\varphi,\boldsymbol\psi,\boldsymbol\varphi',\boldsymbol\psi'$ 
that induce the canonical embeddings
$$
\boldsymbol\varphi_*\colon\cO^E_x\to\cO^{E\oplus F}_x,\qquad \boldsymbol\psi_*\colon\cO_x^F\to \cO_x^{E\oplus F},
$$
resp.~projections
$$
\boldsymbol\varphi'_*\colon\cO_x^{E\oplus F}\to\cO_x^E,\qquad \boldsymbol\psi'_*\colon\cO_x^{E\oplus F}\to\cO_x^F,
$$
the claim follows.
\enddemo

In the categorical language touched upon in 3.2, an ideal structure is an enriched subfunctor of the embedding functor $\bP\to\bO$.

\proclaim{11.3.\ Proposition}(a) For any ideal structure on $\Omega$ and pair of Banach spaces $E,F\neq (0)$, we have $\text{supp } \cO^E/\cJ^E=\text{supp }\cO^F/\cJ^F$.

(b) If the ideal structure is cohesive then $\text{supp }\cO^\bC/\cJ^\bC\subset\Omega$ is closed.
\endproclaim

\demo{Proof}(a)\ Suppose $x\in\Omega$ is not in $\supp\ \cO^E/\cJ^E$.
Any germ $\bold f\in\cO_x^F$ is of form $\boldsymbol\varphi_*\bold e$ with some $\boldsymbol\varphi\in\cO_x^{\Hom(E,F)}$ and $\bold e\in\cO_x^E=\cJ_x^E$, whence $\bold f\in\boldsymbol\varphi_*\cJ^E_x\subset\cJ_x^F$, and so
$x\notin\supp\ \cO^F/\cJ^F$.
Reversing the roles of $E$ and $F$, the claim follows.

(b) If $\ldots\to\cO^E|U@>\varphi>> \cJ^\bC|U\to 0$ is a complete resolution 
over an open $U\subset\Omega$, with $\varphi$ induced by a holomorphic
$f:U\to E^*$, then $U\cap\supp\ \cO^{\bC}/\cJ^{\bC}=\{x\in U\colon f(x)=0\}$ is closed in $U$, whence the claim.
\enddemo

\definition{11.4.\ Definition}By a (complex analytic) subvariety $S$ of $\Omega$ we mean a closed subset $|S|\subset\Omega$ and the specification, for each Banach space $E$, of a sheaf $\cO_S^E$ over $|S|$, so that with some cohesive ideal structure $E\mapsto\cJ^E$
$$
|S|=\supp\ \cO^\bC/\cJ^\bC\qquad\text{ and }\qquad\cO_S^E=(\cO^E/\cJ^E)
\bigr| |S|.
$$
\enddefinition

 From $\cO_\Omega$ the sheaf $\cO_S=\cO_S^\bC$ inherits the structure of a sheaf of rings, and all other sheaves $\cO_S^E$ are modules over it.
The subvariety $S$ uniquely determines the sheaves $\cJ^E$, since $\cO_x^E/\cJ_x^E=(\cO_S^E)_x$ when $x\in |S|$ and $\cJ_x^E=\cO_x^E$ otherwise.

If $\Omega'\subset\Omega$ is open, we denote by $S\cap\Omega'$ the subvariety of $\Omega'$ defined by the ideal structure $\cJ^E|\Omega'$.

We fix a subvariety $S$ of $\Omega$ and the corresponding ideal structure $\cJ^E$.

\subhead 11.5\endsubhead
The sheaves $\cO_S^E$ will be called plain sheaves over $S$.
If $U\subset\Omega$ is open, any section of $\cO_S^{\Hom(E,F)}$ over 
$U\cap |S|$ extends by zero to a section of $\cO^{\Hom(E,F)}/\cJ^{\Hom(E,F)}$ 
over $U$, and it follows from 11.2 that it induces a homomorphism 
$\cO^E/\cJ^E\to\cO^F/\cJ^F$ over $U$, hence a homomorphism $\cO_S^E\to\cO_S^F$ 
over $U\cap |S|$.
Such homomorphisms will be called plain homomorphisms. Germs of plain
homomorphisms form a sheaf of $\cO_S$--modules
$\bHom_{\text{plain}}(\cO_S^E,\cO_S^F)\subset\bHom_\cO(\cO_S^E,\cO_S^F)$;
this sheaf is an epimorphic image of $\cO_S^{\text{Hom}(E,F)}$.

\definition{11.6.\ Definition}An analytic structure on a sheaf $\cA$ of $\cO_S$--modules is the choice, for each Banach space $E$, of a submodule
$\bHom(\cO_S^E,\cA)\subset\bHom_{\cO_S} (\cO_S^E,\cA)$
such that 
\itemitem{(i)}if $\boldsymbol\varphi$ is the germ of a plain homomorphism $\cO_S^E\to\cO_S^F$ at $x$ then
$$
\boldsymbol\varphi^*\bHom(\cO_S^F,\cA)_x\subset\bHom(\cO_S^E,\cA)_x;
$$
\itemitem{(ii)}$\bHom(\cO_S,\cA)=\bHom_{\cO_S}(\cO_S,\cA)$.
\enddefinition

A sheaf of $\cO_S$--modules over $|S|$, endowed with an analytic structure, 
will be called an analytic sheaf over $S$.
A homomorphism $\cA\to\cB$ of analytic sheaves is analytic if the induced maps send each $\bHom(\cO_S^E,\cA)$ in $\bHom(\cO_S^E,\cB)$.
We endow the plain sheaves $\cO_S^F$ with the analytic structure 
$\bHom(\cO_S^E,\cO_S^F)=\bHom_{\text{plain}} (\cO_S^E,\cO_S^F)$.
Thus, for an analytic sheaf $\cA$ over $S$ the sheaf $\bHom(\cO_S^E,\cA)$ 
consists of germs of analytic homomorphisms.

\subhead 11.7\endsubhead
If $\cA$ is any sheaf of $\cO_S$--modules, let $\hat\cA$ denote its extension 
to $\Omega$ by zero outside $|S|$.
This extension clearly has the structure of an $\cO$--module that satisfies 
$\cJ^\bC\hat\cA=0$.
It follows that any $\boldsymbol\psi\in\bHom_{\cO}(\cO,\hat\cA)$ sends 
$\cJ^\bC$ to 0, hence factors through the projection $\cO\to\cO/\cJ^\bC$.

\proclaim{Proposition}(a) If $\cA$ is an analytic sheaf over $S$, then 
$\hat\cA$ has a unique analytic structure with the property that for 
$x\in |S|$ a germ $\boldsymbol\psi\in\bHom_{\cO}(\cO^E,\hat\cA)_x$ is 
analytic if and only if it is the composition of the projection 
$\cO_x^E\to (\cO_S^E)_x$ with a $\boldsymbol\varphi\in\bHom (\cO_S^E,\cA)_x$.

(b) If $\hat\cA$ is endowed with this analytic structure and $U\subset |S|$,
then any analytic homomorphism $\psi:\cO^E|U\to\hat\cA|U$ is the composition
of the projection $\cO^E|U\to\cO^E_S|U$ with a unique analytic homomorphism
$\varphi:\cO_S^E|U\to\cA|U$.
\endproclaim

\demo{Proof}(a) Uniqueness is obvious since $\hat\cA|\Omega\backslash |S|=0$.
The proof of existence consists of verifying that the collection of germs $\boldsymbol\psi$ obtained in the way described above satisfies (i) and (ii) in 3.1, which is straightforward.

(b) Since the projection $\cO^E_x\to(\cO^E_S)_x$ is surjective,
$\boldsymbol\varphi$ in (a) is uniquely determined by $\boldsymbol\psi$,
from which the claim follows.
\enddemo

In what follows, we shall always endow canonical extensions of analytic sheaves over $S$ with the analytic structure described in the above Proposition.
If $\varphi\colon\cA\to\cB$ is an analytic homomorphism of analytic sheaves over $S$ then its extension $\hat\varphi\colon\hat\cA\to\hat\cB$ is also analytic.

\subhead 11.8\endsubhead
We shall call a subset $U\subset |S|$ pseudoconvex if $U=|S|\cap\hat U$ with some pseudoconvex $\hat U\subset\Omega$.

\definition{Definition}A sequence $\cA_\bullet$ of analytic sheaves and homomorphisms over $S$ is completely exact if for each pseudoconvex 
$U\subset |S|$ and for each plain sheaf $\cE$ over $S$ the induced sequence 
$\Gamma(U,\bHom(\cE,\cA_\bullet))$ is exact.
\enddefinition

The following is immediate from 11.7:
\proclaim{Proposition} If $\cA@>\beta>>\cB@>\gamma>>\cC$ is a completely
exact sequence over $S$, then its extension
$\hat\cA@>\hat\beta>>\hat\cB@>\hat\gamma>>\cC$ is completely exact over
$\Omega$.\endproclaim

\definition{Definition}Let $\cA$ be an analytic sheaf over $S$.

(a)\ A complete resolution of $\cA$ is a completely exact sequence
$$
\ldots\to \cF_2\to\cF_1\to\cA\to 0
$$
of analytic homomorphisms, with each $\cF_n$ plain.

(b)\ $\cA$ is cohesive if every $x\in |S|$ has a neighborhood over which $\cA$ admits a complete resolution.
\enddefinition

\subhead 11.9\endsubhead
In the next two theorems we assume $\Omega$ is pseudoconvex.

\proclaim{Theorem~1}If $\ldots @>>> \cA_2 @>>> \cA_1 @>>> \cA_0 @>>> 0$ is 
a completely exact sequence of analytic homomorphisms over $S$ and $\cA_n$ 
has a complete resolution for $n\geq 1$, then so does $\cA_0$.
\endproclaim

\proclaim{Theorem~2}If $0@>>> \cA@>>>\cB @>>>\cC @>>> 0$ is a completely 
exact sequence of analytic homomorphisms over $S$, and two among $\cA,\cB$,
and $\cC$ 
have a complete resolution, then so does the third.
\endproclaim

The proofs are the same as in 4.4 and 4.5.

\proclaim{11.10.\ Theorem}Suppose a Banach space $X$ has a Schauder basis and 
the Hypothesis in 7.1 holds.
Let $S$ be a subvariety of an open $\Omega\subset X$.
If $\cA$ is a cohesive sheaf over $S$, then its canonical extension $\hat\cA$ 
to $\Omega$ is also cohesive.
\endproclaim

\demo{Proof}Cohesion being a local property, we can assume that $\Omega$ is pseudoconvex and $\cA$ has a complete resolution $\cF_\bullet\to\cA\to 0$.
Canonical extension gives rise to a completely exact sequence $\hat\cF_\bullet\to\hat\cA\to 0$.
Since the extension of $\cO_S^F$ is $\cO^F/\cJ^F$, each $\hat\cF_n$ is cohesive by 9.2, and has a complete resolution by 9.1.
Therefore $\hat\cA$ is cohesive by 4.4.
\enddemo

\proclaim{11.11.\ Theorem}Suppose a Banach space $X$ has a Schauder basis and the Hypothesis in 7.1 holds.
Let $\Omega\subset X$ be pseudoconvex and $S$ a subvariety of $\Omega$.
If $\cA$ is a cohesive sheaf over $S$, then
\itemitem{(a)}$\cA$ has a complete resolution; and
\itemitem{(b)}$H^q(|S|,\cA)=0\text{ for }q\geq 1$.
\endproclaim

\demo{Proof}(a)\ By 11.10 the extension $\hat\cA$ is cohesive, so that 9.1 
implies there is a completely exact sequence $\cO^F@>\psi>>\hat\cA\to 0$.
By the Proposition in 11.7 the restriction $\psi\big| |S|$ is the composition of the projection $\cO^F\big| |S|\to\cO_S^F$ with an analytic homomorphism $\varphi\colon\cO_S^F\to\cA$.
We claim that the sequence
$$
\cO_S^F@>\varphi>>\cA\to 0\tag11.2
$$
is completely exact.

Indeed, let $\hat U\subset\Omega$ be pseudoconvex, $U=\hat U\cap |S|$, and $E$ a Banach space.
Any $\alpha\in\Gamma (U,\bHom(\cO_S^E,\cA))$ extends to an $\hat\alpha\in\Gamma(\hat U,\bHom(\cO^E/\cJ^E,\hat\cA))$, which, composed with the projection $\cO^E\to\cO^E/\cJ^E$ gives a $\beta\in\Gamma (\hat U,\bHom(\cO^E,\hat\cA))$, again by 11.7.
Since $\cO^F @>\psi>>\hat\cA\to 0$ was completely exact, $\beta=\psi_*\gamma$ with some $\gamma\in\Gamma (\hat U,\bHom (\cO^E,\cO^F))$.
As $\gamma$ maps $\cJ^E|\hat U$ in $\cJ^F|\hat U$, it descends to a section of $\bHom(\cO^E/\cJ^E,\cO^F/\cJ^F)$ over $\hat U$, whose restriction to $U$,
$$
\delta\in\Gamma (U,\bHom (\cO_S^E,\cO_S^F)),
$$
satisfies $\varphi_*\delta=\alpha$.
This proves (11.2) is completely exact, as claimed. 

Next, Theorem~2 in 11.9 implies that $\Ker\varphi$ in the completely exact sequence
$$
0\to\Ker\varphi\to\cO_S^F@>\varphi>> \cA\to 0
$$
is cohesive.
Repeating the above construction with $\Ker\varphi$ instead of $\cA$ we obtain a completely exact sequence $0\to\Ker\varphi'\to\cO_S^{F'} @>\varphi'>>\Ker\varphi\to 0$, with $\Ker\varphi'$ cohesive, and so on.
Consolidating the short sequences obtained in this way, we construct a complete resolution
$$
\ldots\to\cO_S^{F'}\to \cO_S^F\to\cA\to 0.
$$
(b)\ The extension $\hat\cA$ being cohesive, by 9.1
$$
H^q (|S|,\cA)\approx H^q(\Omega,\hat\cA)=0\qquad\text{ for }q\geq 1.
$$
\enddemo

\head 12.\ Two Applications\endhead
\subhead 12.1\endsubhead
As explained in 5.2, cohesive sheaves pulled back by biholomorphisms stay cohesive.
For this reason one can define cohesive (and plain and analytic) sheaves over an arbitrary rectifiable complex Banach manifold.

Now consider a Banach space $X$ that has a Schauder basis and satisfies the Hypothesis in 7.1.
Let $\Omega\subset X$ be open.
A direct submanifold $M\subset\Omega$ determines an ideal structure, with $\cJ^E\subset\cO^E=\cO^E_\Omega$ consisting of germs that vanish on $M$.
By 10.3 this structure is cohesive, and so it defines a subvariety 
$S$ with support $|S|=M$.
On $M$ we have two notions of plain sheaves:\ the sheaves $\cO_M^E$ of $E$ valued holomorphic germs on $M$, and the sheaves $\cO_S^E=(\cO^E/\cJ^E)|M$ discussed in Section~11.
Similarly, there are two notions of plain homomorphisms.
However, restricting germs in $\cO^E$ to $M$ induces an isomorphism $\cO_S^E\approx\cO_M^E$, which isomorphism intertwines plain homomorphisms $\cO_S^E\to\cO_S^F$ and $\cO_M^E\to\cO_M^F$.
It follows that any sheaf of $\cO_S$--modules has a canonical structure of an $\cO_M$--module, and vice versa; any analytic sheaf over $S$ has a canonical structure of an analytic sheaf over $M$, and vice versa; analytic homomorphisms will be the same, whether the analytic sheaves are considered over $S$ or $M$; and finally, whether a sheaf is cohesive does not depend on whether it is considered over $S$ or over $M$.

\subhead 12.2\endsubhead
The following theorem is therefore a special case of 11.11:

\proclaim{Theorem}Suppose a Banach space $X$ has a Schauder basis and the 
Hypothesis in 7.1 holds, $\Omega\subset X$ is pseudoconvex, and 
$M\subset\Omega$ is a direct submanifold.
If $\cA$ is a cohesive sheaf over $M$ then
\itemitem{(a)}$\cA$ has a complete resolution;
\itemitem{(b)}$H^q(M,\cA)=0$ for $q\geq 1$.
\endproclaim

\proclaim{12.3.\ Theorem}Let $X,\Omega$, and $M$ be as above, and $\nu=(T\Omega|M)/TM$ the normal bundle.
Any neighborhood of $M$ contains a pseudoconvex neighborhood $O\subset\Omega$

(a)\ which holomorphically retracts on $M$; and

(b)\ there is a biholomorphism between $O$ and a neighborhood of $M\subset\nu$ ($M$ embedded as the zero section in $\nu$), that is the identity on $M$.
\endproclaim

This generalizes the Docquier--Grauert theorem, a special case of a 
theorem of Siu, see [DG,Su], and [L6, Theorem~1.4].
In the former two $\dim X<\infty$; in the latter $M$ was assumed 
biholomorphic to an open subset of a Banach space. Particular cases of
both 12.2 (b) and 12.3 were also proved by Patyi, see [P5, Theorem 6.2].

\demo{Proof}The key is the vanishing $H^1(M,\Hom(\nu,TM))=0$, which follows from 12.2 if one takes $\cA$ to be the sheaf of holomorphic sections of the locally trivial holomorphic Banach bundle $\Hom(\nu,TM)$.
This, combined with the exact sequence of holomorphic Banach bundles
$$
0\to\Hom(\nu,TM)\to\Hom(\nu,T\Omega|M)\to\Hom(\nu,\nu)\to 0
$$
provides a subbundle $\nu'\subset T\Omega|M$, complementary to $TM$, as 
in [DG, L6]. From here $O$ and its biholomorphic embedding in $\nu$ are 
constructed as in [L6], to which we refer the reader for details.
\enddemo

\head 13.\ Appendix\endhead
\subhead 13.1\endsubhead
If $\Omega$ is an open subset of a finite dimensional Banach space $X$ and $S\subset\Omega$ is an analytic subset, then the sheaf $\cJ\subset\cO_\Omega$ of germs vanishing on $S$, the ideal sheaf of $S$, is coherent.
Here we shall address the corresponding problem in infinite dimensional $X$.
Let $\Omega\subset X$ be open.
A subset $S\subset\Omega$ is called analytic if $\Omega$ can be covered by 
open sets $U\subset X$, and for each $U$ there are a Banach space $E$ and 
a holomorphic function $f\colon U\to E$ such that
$$
S\cap U=\{x\in U\colon f(x)=0\},
$$
see [D1-2,R].
Analytic sets can be pathological:\ for example, Douady pointed out
in [D2] that any compact metric space can be homeomorphically embedded in 
some Banach space as an analytic subset.
A variant of Douady's construction shows that any compact differential 
manifold, possibly with boundary, can be homeomorphically embedded in a 
Hilbert space as an analytic subset.
The notion of an analytic set is clearly too generous, and in complex 
geometry one should restrict to a smaller class of sets.
We propose that the correct class will consist of those sets 
$S\subset\Omega$ whose ideal sheaf is cohesive.
We shall show below that this is a genuine restriction:\ even in spaces with 
unconditional basis there are analytic subsets whose ideal sheaf is not 
cohesive.

\subhead 13.2\endsubhead
Let $X$ denote one of the spaces $l^p, 1\leq p<\infty$, or $c_0$, and consider 
the following embedding $\alpha$ of the closed unit disc $\Delta\subset\bC$ 
into $X$
$$
\alpha(\zeta)=(\zeta^{n!}/n^2)^\infty_{n=1}\in X,\qquad \zeta\in\Delta.
$$
The image $\alpha(\Delta)\subset X$ is analytic, since it is the zero set of 
the holomorphic map
$$
X\ni (x_n)\mapsto \bigl((n^2x_n-x_1^{n!})/n^{n!}\bigr)\in X.
$$
Let $\cJ\subset\cO_X=\cO$ denote the sheaf of germs vanishing on 
$\alpha(\Delta)$.

\proclaim{Theorem}For no Banach space $F$ and analytic homomorphism 
$\cO^F\to\cJ$ is the restriction $\cO_0^F\to\cJ_0$ surjective.
\endproclaim

Therefore $\cJ$ is not cohesive because it cannot be included in an exact 
sequence of analytic homomorphisms $\cO^F\to\cJ\to 0$, cf.~9.1.---To 
prepare the proof, consider the following continuous action of the circle 
$S^1=\bR/2\pi\bZ$ on $X$
$$
g_t x=(e^{in!t}x_n),\qquad t\in S^1,\ x=(x_n)\in X.\tag13.1
$$
If $E$ is a Banach space, any holomorphic function $h\colon B(R)\to E$ has
a Fourier series
$$
h\sim\sum^\infty_{k=0} h_k,\quad h_k=
{1\over 2\pi}\int_0^{2\pi} e^{-ikt} g_t^* h dt.\tag13.2
$$
The terms $h_k$ are holomorphic on $B(R)$ and satisfy $g_t^* h_k=e^{ikt}h_k$.
If
$$
\gather
h(x_1,\dots,x_m,0,0,\dots)=
\sum_{j_1,\dots,j_m\ge 0} a_{j_1\ldots j_m} x_1^{j_1}\cdots x_m^{j_m},
\qquad\text{then}\\
h_k (x_1,\dots,x_m,0,0,\dots)=
\sum_{{\sum \mu!j_\mu}=k} a_{j_1\ldots j_m} x_1^{j_1}\cdots x_m^{j_m}.
\endgather
$$
The last expression is a polynomial, independent of the variables 
$x_\mu$ with $\mu!>k$.
It follows that $h_k(x_1,x_2,\ldots)$ itself is a polynomial, depending 
only on the variables $x_\mu$ with $\mu!\leq k$.
It also follows that $h(x)=\sum_k h_k(x)$ if the sequence $x\in B(R)$ has 
only finitely many nonzero terms.

\demo{Proof of the Theorem} Consider an analytic homomorphism
$\cO^F\to\cJ$ induced by a holomorphic 
$\varphi\colon X\to F^*$ that vanishes on $\alpha(\Delta)$.
Let $\boldsymbol\varphi$ be the germ at $0$ of $\varphi$ and
$\bold u\in\cJ_0$ be the germ at 0 of
$$
u(x)=\sum^\infty_{n=1}(n^{-4}x_1^{n!}-n^{-2}x_n),\qquad x\in B(1).
$$
We shall show that there is no $\bold f\in\cO_0^F$ solving 
$\boldsymbol\varphi\bold f=\bold u$.

Suppose there were.
The Fourier components of $\boldsymbol\varphi,\bold f,\bold u$ would then 
satisfy
$$
\sum^{n!}_{k=0}\boldsymbol\varphi_k(x) \bold f_{n!-k}(x)
=\bold u_{n!}(x),\qquad n\in\bN.\tag13.3
$$
Since $\alpha(\Delta)$ is $g_t$ invariant, all $\boldsymbol\varphi_k$
vanish on $\alpha(\Delta)$. We substitute
$$
x=\alpha^{(n)}(\zeta)=
\biggl(\zeta,\frac{\zeta^{2!}}{2^2},\frac{\zeta^{3!}}{3^2},\dots,
\frac{\zeta^{(n-1)!}}{(n-1)^2},0,0,\dots\biggr).
$$
As $\boldsymbol\varphi_k(x)$ depends only on $x_1,\dots,x_{n-1}$ if
$k<n!$, we have $\boldsymbol\varphi_k(\alpha^{(n)}(\zeta))=
\boldsymbol\varphi_k(\alpha(\zeta))=0$. By (13.2)
$\bold f_0$ is the constant $\bold f(0)$. Therefore we
obtain from (13.3)
$$
\boldsymbol\varphi_{n!}\bigl(\alpha^{(n)}(\zeta)\bigr)\bold f(0)=
\bold u_{n!}\bigl(\alpha^{(n)}(\zeta)\bigr).
$$
Both sides are germs of polynomials at $0\in\bC$. If they agree,
the corresponding polynomials agree on all of $\bC$. In particular,
since $u_{n!}(x)=n^{-4}x_1^{n!}-n^{-2}x_n$,
$$
\varphi_{n!}\bigl(\alpha^{(n)}(n^{5/n!})\bigr)\bold f(0)=n.\tag13.4
$$

By the dominated convergence theorem
$\lim_n\alpha^{(n)}(n^{5/n!})=(1/m^2)_m$. Thus the set
$\{\alpha^{(n)}(n^{5/n!})\colon n\in \Bbb N\}$ has compact closure, and
so does its $g_t$ orbit. If $K$ is the closure of the orbit, then from 
(13.2)
$$
||\varphi_{n!}\bigl(\alpha^{(n)}(n^{5/n!})\bigr)||_{F^*}\le\max_K||\varphi||_{F^*},
$$
hence by (13.4) $||\bold f(0)||_F\max_K||\varphi||_{F^*}\ge n$.
This, however, cannot hold for all $n$, which proves that the equation 
$\boldsymbol\varphi\bold f=\bold u$ has no solution $\bold f$.
\enddemo

\subhead 13.3\endsubhead
More regular examples of analytic sets can exhibit the same phenomenon.
Consider the following $C^\infty$ embedding
$$
\beta\colon\Delta\ni\zeta\mapsto (\zeta^{n!}/(n!)^{\log n})\in X.
$$
The image $\beta(\Delta)$ is an analytic subset of $X$; but a variant of the 
above reasoning shows that its ideal sheaf is not cohesive, either.

\enddocument
\bye

\bye